\documentclass[11pt]{amsart}
\usepackage[utf8]{inputenc}
\usepackage{amsmath,amssymb,amsthm, mathrsfs, mathtools}
\usepackage{amsfonts,amscd,verbatim,calc,amsgen}
\usepackage{indentfirst}
\usepackage{comment,relsize,braket}
\usepackage{enumerate}
\usepackage[numbers]{natbib}
\usepackage{graphicx}
\usepackage{mathrsfs}
 \usepackage{textcomp}
\usepackage{lscape}
\usepackage{mathpazo}
\usepackage{color}

\usepackage[all]{xy}
\usepackage{fullpage}
\usepackage{tikz}
\usepackage[top=1in, bottom=1in, left=1in, right=1in]{geometry}

\usepackage{fancyhdr}
\DeclareMathOperator{\HH}{H}

\DeclareMathOperator{\hgt}{ht}

\DeclareMathOperator{\im}{Im}

\DeclareMathOperator{\Soc}{Soc}

\DeclareMathOperator{\depth}{depth}

\DeclareMathOperator{\chr}{char}

\newcommand{\ncom}{\newcommand}
\ncom{\ov}{\overline}
\ncom{\N}{\mathbb N}
\ncom{\Q}{\mathbb Q}
\ncom{\Z}{\mathbb Z}
\ncom{\F}{\mathbb{F}} 
\ncom{\C}{\mathbb{C}} 
\ncom{\p}{\mathfrak{p}} 
\ncom{\q}{\mathfrak{q}} 
\ncom{\m}{\mathfrak{m}}
\ncom{\n}{\mathfrak{n}}
\ncom{\kk}{\mathbf{k}}
\ncom{\mL}{\mathcal{L}} 
\ncom{\mN}{\mathcal{N}} 
\ncom{\f}{\mathcal{F}}
\ncom{\g}{\mathcal{G}}
\ncom{\I}{\mathcal{I}}
\ncom{\J}{\mathcal{J}} 
\ncom{\R}{\mathcal{R}}
\ncom{\Sa}{\mathcal{S}}
\ncom{\A}{\mathcal{A}}
\ncom{\B}{\mathcal{B}}
\ncom{\cc}{\mathbf{C}}
\ncom{\hc}{\mathbf{HC}}

\def\thebibliography#1{\section*{References}
\list{[\arabic{enumi}]}{\settowidth \labelwidth{[#1]} \leftmargin
\labelwidth \advance \leftmargin \labelsep \usecounter{enumi}}
\def\newblock{\hskip .11em plus .33em minus .07em} \sloppy
\clubpenalty 4000 \widowpenalty 4000 \sfcode`\.=1000 \relax}

\newtheorem{Theorem}{Theorem}[section]
\newtheorem{Lemma}[Theorem]{Lemma}
\newtheorem{Corollary}[Theorem]{Corollary}
\newtheorem{Proposition}[Theorem]{Proposition}

\newtheorem{Remarks}[Theorem]{Remarks}
\newtheorem{Example}[Theorem]{Example}

\newtheorem{Definition}[Theorem]{Definition}
\newtheorem{Notation}[Theorem]{Notation}

\newtheorem{Conjecture}[Theorem]{Conjecture}

\numberwithin{equation}{section}

\begin{document}
\title[]{On the Vanishing of \\ [2mm] the normal  Hilbert coefficients  of ideals}
\thanks{The first author is supported by a UGC fellowship, Govt. of India}

\author[]{ Kriti Goel}
\author[]{Vivek Mukundan}
\author[]{J. K. Verma}
\address{Indian Institute of Technology Bombay, Mumbai, INDIA 400076}
\address{University of Virginia, Charlottesville, VA 22904, USA}
\email{kriti@math.iitb.ac.in}
\email{vm6y@virginia.edu}
\email{jkv@math.iitb.ac.in}

\begin{abstract}
	Using vanishing of graded components of local cohomology modules of the Rees algebra of the normal filtration of an ideal, we give bounds on the normal reduction number. This helps to get necessary and sufficient conditions in Cohen-Macaulay local rings of dimension $d\geq 3$, for the vanishing of the normal Hilbert coefficients $\ov{e}_k(I)$ for $k\leq d,$ in terms of the normal reduction number. 
\end{abstract}
\maketitle

\section{Introduction}
The objective of this paper is to explore vanishing of normal Hilbert coefficients and bounds on normal reduction numbers, inspired by Itoh's conjecture \cite[Section 5]{itoh92}. In order to describe the theme, we recall the notions of normal Hilbert polynomial and normal reduction number.

Let $I$ be an ideal of a commutative ring $R.$ We say that $x\in R$ is  integral over $I$ if it satisfies an equation
\[x^n+a_1x^{n-1}+\dots+a_n=0\]
for some $a_i\in I^i$, for all $i=1,2,\dots,n.$ The integral closure of $I,$ denoted by $\ov{I},$ is the $R$-ideal consisting of all $x\in R$ which are integral over $I.$  Let $\N$ denote the set of non-negative integers. The filtration of ideals $\{\ov{I^n}\}_{n\in \N}$ is called the normal filtration of $I.$ David Rees \cite{reesAUrings} showed that if $(R,\mathfrak{m})$ is an analytically unramified (i.e., the $\mathfrak{m}$-adic completion $\widehat{R}$ is reduced) local ring, then there exists an $h\in \N$ such that for all $n\geq 0,$ 
\begin{eqnarray}\label{reesadm}
\ov{I^{n+h}}\subset I^n.
\end{eqnarray} 
Let  $I$ be an $\m$-primary ideal of a Noetherian local ring $(R,\m).$ The function $\ov{H}_I(n)=\ell(R/\ov{I^{n}})$ is called the {\em normal Hilbert function of $I.$} 
Using (\ref{reesadm}), Rees showed that if $R$ is an analytically unramified Noetherian local ring of dimension $d$, then there exists a polynomial $\ov{P}_I(x)\in \Q[x]$ of degree $d$ such that for all large $n,$ $\ov{H}_I(n)=\ov{P}_I(n).$ The polynomial $\ov{P}_I(x)$ is called the {\em normal Hilbert polynomial of $I.$} We write it as

\[ \ov{P}_I(x) = \ov{e}_0(I) \binom{x+d-1}{d} - \ov{e}_1(I) \binom{x+d-2}{d-1} + \cdots + (-1)^d \ov{e}_d(I). \]
The coefficients $\ov{e}_0(I),\dots,\ov{e}_d(I)$ are integers and are called the {\em normal Hilbert coefficients of $I.$} Let $R/\m$ be infinite. Then there exists a minimal reduction  $J$  of $I$ and $r\in \N$ such that $J\ov{I^n}=\ov{I^{n+1}}$ for all $n\geq r.$ The minimal such $r,$ denoted by $\ov{r}_J(I)$, is called the normal reduction number of $I$ with respect to $J.$ If $I$ is generated by a system of parameters, then we shall write $\ov{r}(I)=\ov{r}_I(I).$ One of the recurring themes in the investigation of the normal Hilbert coefficients has been their relationship with the normal reduction number of the normal filtration of $I.$ The following conjecture of Itoh is open since 1992:

\begin{Conjecture}[\bf Shiroh Itoh \cite{itoh92}] \label{introconj}
	Let $R$ be a Gorenstein local ring of dimension $d\geq 3.$ Let $I$ be an $\m$-primary ideal. Then
	\[\ov{e}_3(I)=0 \ \text{ if and only if } \ \ov{r}(I)\leq 2.\]
\end{Conjecture}

Itoh settled this conjecture in the affirmative for the normal filtration of $\m.$ See the papers by Corso-Polini-Rossi \cite{cmr} and Kummini-Masuti \cite{kuma} for recent progress on this conjecture. The following result of Huneke and Itoh, which we refer to as the {\em Huneke-Itoh Intersection Theorem}, plays an important role in the study of the normal Hilbert polynomial.

\begin{Theorem}[\bf Huneke \cite{huneke1987}, Itoh \cite{itoh88}] \label{introHunekeItoh}
	Let $R$ be a Noetherian ring and $I$ be generated by a regular sequence. Then for all $n\geq 0,$
	\[I^n \cap \ov{I^{n+1}}=\ov{I}I^n. \]
\end{Theorem}
Itoh, very effectively, exploited the vanishing of graded components of the local cohomology module $\HH^2_J(\ov{\R'}(I))$ to prove the above result (also see \cite{hongUlrich}). Here $\ov{\R'}(I)$ represents the integral closure of the extended Rees algebra $\R'(I)$ and $J=(t^{-1},It)$ is an $\R'(I)$-ideal. In fact, the strength of the above theorem lies in the fact that, along with $\ov{r}(I)\leq 2$, it automatically implies that $\ov{e}_3(I)=0$ \cite[Proposition 10]{itoh88}. Subsequently, in \cite{itoh92}, he used the machinery of general extensions to show the bounds for $\ov{e}_1(I),\ov{e}_2(I)$. Further he showed that these bounds are achieved when the reduction number is bounded. Since then, there have been many results on the nature of $\ov{e}_1(I),\ov{e}_2(I)$ (\cite{Huckaba1,huckabaMarley,huneke1987,mandal2012chern}). So it is natural to ask if similar behavior is exhibited for higher Hilbert coefficients. Such characterizations seem possible when the associated graded ring has depth at least $\dim R-1$ \cite{Huckaba1,huckabaMarley,marley89,marleyThesis}.

Inspired by the above evidence on the relationship between the graded components of the local cohomology modules and the normal Hilbert coefficients, an analogue of Conjecture \ref{introconj} for the higher normal Hilbert coefficients $\ov{e}_k(I), k\geq 3$ can be asked. In order to prove such an analogue we  use \emph{the condition $HI_r$} introduced in \cite{GMV}. Recall that
\begin{Definition}\label{intronormalHIdef}{\rm
		For $r \ge 1,$ the normal filtration $\f=\{\ov{I^n} \}$ is said to satisfy \emph{the condition $HI_{r}$} if 
		\begin{align*}
		I^n \cap \ov{I^{n+r}} = I^n \ov{I^r}\text{ for all }n \geq 0.
		\end{align*}
}\end{Definition}
Notice that \emph{the condition $HI_1$} is the same as that of the equality in Theorem \ref{introHunekeItoh}. Interestingly, as it happens in the case of Theorem \ref{introHunekeItoh}, the normal filtration $\f=\{ \ov{I^n}\}$ satifies \emph{the condition $HI_r$} when some graded components of the local cohomology modules $\HH^i_J(\ov{\R'}(I))$ vanish. We prove the following result in section 3. 

\begin{Theorem}\label{introGITintclfilt}
	Let $R$ be an equidimensional, universally catenary, and an analytically unramified Noetherian local ring of dimension $d.$ Let $I$ be an ideal generated by an $R$-regular sequence. Then \\
	{\rm(1)} If $\hgt(I)=1,$ then $\f$ satisfies  $HI_r$  for all $r \ge 1.$ \\
	{\rm(2)} Let $\hgt(I) \ge 2$ and for some $r \geq 1$, let $\HH^i_J(\ov{\R'}(I))_{j} = 0$ for all $i,j$ such that $i+j=r+1$ and $0 \le i \le \hgt(I).$ Then $\f$ satisfies {\em the condition $HI_{r}$}.
\end{Theorem}

To prove Theorem \ref{introGITintclfilt}, we use the method of specialization of integral closure used by J. Hong and B. Ulrich in \cite{hongUlrich} (refer to Theorem \ref{specialprime} for a general version). In fact, we show as an evidence, Example \ref{exampleHunekeHuckaba}, where the non-vanishing of the graded pieces of the local cohomology modules appears in conjunction with the violation of \emph{the condition $HI_2$}.

In the interluding section 4, we explore the relationship between vanishing of local cohomology modules of the extended Rees algebra of the normal filtration and the normal reduction number. The aim of the section is to prove the following theorem.

\begin{Theorem}
	Let $(R,\m)$ be a $d$-dimensional analytically unramified, Cohen Macaulay local ring with $d\geq 2$ and let $I$ be a parameter ideal. Suppose that $\HH^i_{J}(\ov{\R'}(I))_j=0$ for all $i,j$ such that $3 \leq i+j \leq k-1$ and $2 \leq i \leq d.$ If  
	\[ \ov{e}_2(I)=(k-2)\ov{e}_1(I)-\sum_{i=0}^{k-3}(k-2-i) \ \ell\left(\frac{\ov{I^{i+1}}}{I\ov{I^i}}\right), \] 
	then $\ov{r}(I)\leq k-1$.
\end{Theorem}
The above result can be viewed as a generalization of Itoh's characterization \cite[Theorem 2]{itoh92} for ideals of reduction number at most two.

With the required tools in place, we prove the analogue of Conjecture of \ref{introconj} for higher Hilbert coefficients in section 5. For this purpose, we start by generalizing Itoh's result \cite[Proposition 10]{itoh88}, for normal Hilbert polynomial of parameter ideals with normal reduction number at most $k-1.$

\begin{Proposition}
	Let $R$ be a $d$-dimensional analytically unramified, Cohen-Macaulay local ring and $I$ be a parameter ideal. Let $\f=\{ \ov{I^n} \}$ satisfy {\em the conditions $HI_p$} for all $p \leq k-2$ and let $k\leq d-1.$ Then for all $n \geq k-2,$
	{\small\begin{equation} \label{poly1}
		\ell\left(\frac{R}{\ov{I^{n+1}}}\right)
		\leq \ell\left(\frac{R}{I}\right)\binom{n+d}{d} - \alpha_1(\f)\binom{n+d-1}{d-1} + \cdots + (-1)^{k-1} \alpha_{k-1}(\f)\binom{n+d-(k-1)}{d-(k-1)}
	\end{equation}}
	where $\displaystyle \alpha_j(\f) = \sum_{i=j-1}^{k-2} \binom{i}{j-1}\ell(\ov{I^{i+1}}/I\ov{I^i})$ for all $j=1,\ldots,k-1.$ The equality holds in the equation {\rm(\ref{poly1})} if and only if $\ov{r}(I) \leq k-1.$ In this case, $\ov{G}(I)$ is Cohen-Macaulay. 
\end{Proposition}	
The above proposition also shows that, under the conditions, $\ov{r}(I)\leq k-1$ implies that the normal Hilbert coefficient $\ov{e}_k(I)$ vanish. The converse is what is sought after. When the associated graded ring is almost Cohen-Macaulay, the positivity of $\ov{e}_k(I)$ is well known \cite[Corollary 2.9]{marleyThesis}. But in Theorem \ref{vanishingOfe_k}, we discuss the positivity of the normal Hilbert coefficient $\ov{e}_k(I)$ under weaker hypothesis and also give the implications of the vanishing of $\ov{e}_k(I)$. Now using the above proposition, we are able to characterize the vanishing of $\ov{e}_k(I)$ in the following theorem. %Notice that this also generalizes the result in \cite{cmr}. 

\begin{Theorem} 
	Let $(R,\m)$ be a $d$-dimensional analytically unramified, Cohen-Macaulay local ring with $d\geq 3$, and let $I$ be a parameter ideal such that $\overline{I}=\m$. Suppose that for some $k \leq d,$ $\HH^i_{J}(\ov{\R'}(I))_j=0$ for all $i,j$ such that $3 \leq i+j \leq k-1$ and $2 \leq i \leq d$ and $\ell(\ov{I^{k-1}}/I\ov{I^{k-2}})\geq t(R)$, where $t(R)$ denotes the type of the ring $R$. Then $\ov{e}_k(I)=0$ if and only if $\ov{r}(I)\leq k-1.$ In this case, $\ov{G}(I)$ is Cohen-Macaulay.
\end{Theorem}

Finally, the effectiveness of \emph{the condition $HI_r$} for normal filtration is reason enough to extend Theorem \ref{introGITintclfilt} for any admissible filtration. Recall that in a ring $R$, a filtration $\f$ of ideals  is a descending chain $R=I_0 \supseteq I_1 \supseteq I_2 \supseteq \cdots$ of ideals such that $I_iI_j \subseteq I_{i+j}$ for all $i,j \geq 0.$ Let $I$ be an ideal of $R.$ Then a sequence of ideals $\f=\{I_n\}_{n \in \Z}$ is called an $I$-filtration if in addition, $I^n \subseteq I_n$ for all $n.$ An $I$-filtration is called $I$-admissible if there exists a $k \in \N$ such that $I_n \subseteq I^{n-k}$, for all $n.$ If $\f=\{I_n\}_{n\in \Z}$ is an $I$-admissible filtration, then the extended Rees algebra $\R'(\f)=\oplus_{n\in \Z}I_nt^n$ and the Rees algebra $\R(\f)=\oplus_{n\in \N}I_n t^n$ are finite modules over the extended Rees algebra $\R'(I)=\oplus_{n \in \N}I^nt^n$ and the Rees algebra $\R(I)=\oplus_{n \in \N}I^nt^n$ respectively. Let $G(\f)=\bigoplus_{n\ge 0}I_n/I_{n+1}$ be the associated graded ring of an $I$-admissible filtration $\f=\{I_n\}.$ 
{\em condition $HI_r$} is then defined as follows:
\begin{Definition}{\rm
		Let $R$ be a Noetherian  ring and $I$ be an $R$-ideal. Let $\f=\{I_n\}_{n \in \Z}$ be an $I$-admissible filtration, where $I_n=R$ for all $n \leq 0.$ For $r \ge 1,$ the filtration $\f$ is said to satisfy \emph{the condition $HI_{r}$} if for all $n \geq 0$, 
		\begin{align*}
			I^n \cap I_{n+r} = I^n I_r.
		\end{align*}
}\end{Definition}

In section 7, we find the conditions under which an admissible $\mathcal{F}$ satisfies {\em the condition $HI_r.$}
\begin{Theorem}
	Let $R$ be a $d$-dimensional Cohen-Macaulay, local ring and let $I$ be an $R$-ideal generated by an $R$-regular sequence of height $g \ge 1.$ Let $\f$ be an $I$-admissible filtration. If $\HH^i_J(\R'(\f))_{j}=0$ for all $i,j$ such that $i+j=r+1$ and $1 \le i \le g$, and $\depth G(\f) \ge g-1$, then $\f$ satisfies {\em the condition $HI_r.$}
\end{Theorem}

\section{Settings, Notation and Preliminaries}
In this section we introduce the settings and notation which we use throughout the article. At the end of the section, for ease of reference, we also present the theorems we often use in the article.

\begin{Notation}\label{notation}\rm{
\begin{enumerate}[(a)]
	\item $(R,\mathfrak{m})$ is a Noetherian local ring, $I=(a_1,\ldots,a_g)$ is an $R$-ideal.
	\item $t$ is an indeterminate, $\mathcal{R}'(I) = \bigoplus_{n \in \Z} I^nt^n$ is the extended Rees algebra of $I$ and $\mathcal{R}(I) = \bigoplus_{n \in \N} I^nt^n$ is the Rees algebra of $I.$ Put $\mathcal{R}(I)_+ = \bigoplus_{n > 0} I^nt^n$ and $J=(t^{-1},It)\R'(I).$
	\item $\mathcal{F}=\{I_n\}$ is an $I$-admissible filtration.
	$\mathcal{R}'(\mathcal{F}) = \bigoplus_{n\in \Z} I_nt^n$ is the extended Rees algebra and $\mathcal{R}(\mathcal{F}) = \bigoplus_{n \in \N} I_nt^n$ is the Rees algebra of $\mathcal{F}.$
	$G(\f)=\bigoplus_{n\ge 0}I_n/I_{n+1}$ is the associated graded ring of $\f.$
	\item If $\f$ is the normal filtration of $I$, i.e., $\f=\{\ov{I^n} \}$, then the extended Rees algebra of $\f$ is denoted by $\ov{\R'}(I)$ and the Rees algebra of $\f$ is denoted by $\ov{\R}(I).$
	\item Let $z_1,\dots,z_g$ be indeterminates over $R.$ Then $R' = R[z_1,\dots,z_g]$ and $I' = IR'.$
	 $\mathcal{R}'(I') = \bigoplus_{n \in \Z} (I')^nt^n$ is the extended Rees algebra, $\mathcal{R}(I') = \bigoplus_{n \in \N} (I')^nt^n$ is the Rees algebra of $I'$ and $\ov{\mathcal{R}'}(I')$ is the extended Rees algebra of the filtration $\{\ov{(I')^n}\}.$
	\item $R'' = R[z_1,\dots,z_g]_{\m R[z_1,\dots,z_g]}$ is the general extension, $I'' = IR''$ and $x = \sum_{i=1}^{g}z_ia_i \in R''.$
	 $\mathcal{R}'(I'') = \bigoplus_{n \in \Z} (I'')^nt^n$ is the extended Rees algebra of $I''$ and $J''=(t^{-1},I''t) \mathcal{R}'(I'')$ is an $\mathcal{R}'(I'')$-ideal. Put $\mathcal{R}(I'')_+ = \bigoplus_{n > 0} (I'')^nt^n.$
	\item $T=R''/(x)$ and $IT=I''/(x).$
\end{enumerate}
}
\end{Notation}
\begin{Definition}\label{Crnotation}
	Let $\f$ be an $I$-admissible filtration. Then for $r\in\mathbb{N}$,
	\begin{align*}
	\f \text{ satisfies }\cc_{r}, \text{ if }\HH^i_J(\mathcal{R}'(\f))_j=0\text{ for }i+j=r+1\text{ and } 0\leq i\leq \hgt(I)
	\end{align*}
	and 
		\begin{align*}
	\f \text{ satisfies }\hc_{r}, \text{ if }\HH^i_{\mathcal{R}(I)_+}(\mathcal{R}'(\f))_j=0\text{ for }i+j=r+1\text{ and } 0\leq i\leq \hgt(I).
	\end{align*}
\end{Definition}
Using the above notation, notice that if $\f=\{\ov{I^n} \}$ is the normal filtration, then $\f$ satisfies $\cc_{r}$ if $\HH^i_J(\ov{\mathcal{R}'}(I))_j=0\text{ for }i+j=r+1$ and $0\leq i\leq \hgt(I)$.

\begin{Theorem}{\rm\cite[Theorem 2.1]{hongUlrich}}
	Let $R$ be a Noetherian, locally equidimensional, universally catenary ring such that $R_{red}$ is locally analytically unramified. Let $I=(a_1,\dots,a_g)$ be an $R$-ideal of height at least $2$, $I'=IR'$ and $x=\sum a_iz_i$. Then  $\overline{I'/(x)}=\overline{I'}/(x)$.
\end{Theorem}

\begin{Theorem}{\rm \cite[Proposition 13]{itoh92}}
	Let $(R,\m)$ be an analytically unramified, Cohen-Macaulay local ring of dimension $d \geq 2$. Let $I=(a_1,\dots,a_d)$ be a parameter ideal. Then for all $i=0,1,\dots,d-1$,
	\[\HH^i_{\R(I)_+}(\ov{\R'}(I)) \simeq \HH^i_J(\ov{\R'}(I))\]
	and there is an exact sequence
	\begin{equation} \label{itohseq}
	0 \rightarrow \HH^d_J(\ov{\R'}(I)) \rightarrow  \HH^d_{\R(I)_+}(\ov{\R'}(I)) \rightarrow \HH^d_{\m}(R)[t,t^{-1}] \rightarrow \HH^{d+1}_J(\ov{\R'}(I)) \rightarrow 0.
	\end{equation}
\end{Theorem}
Thus by using the above proposition, it is clear that, if the normal filtration $\f=\{\ov{I^n} \}$ satisfies $\hc_{r}$ then $\f$ satisfies $\cc_{r}$. The converse is true when $\hgt(I)<d$.

\begin{Lemma}{\rm \cite[Lemma 14]{itoh92}}
Under the conditions of the above theorem, suppose $r$ is an integer such that $\HH^i_{\R(I)_+}(\ov{\R'}(I))_j=0$ for $i+j=r+1$, then $\HH^i_{\R(I)_+}(\ov{\R'}(I))_j=0$ for $i+j\geq r+1$ and $\ov{I^{n+r}}=I^n\ov{I^r}$ for every $n\geq 0$.
\end{Lemma}
%\cite[Theorem 2]{itoh88}

\section{The condition $HI_r$ for normal filtrations}\label{normafiltrationHI}

In this section we prove the first main result, Theorem \ref{integralclosure}, of this article. In here, we discuss the necessary conditions for the normal filtration $\f=\{\ov{I^n} \}$ to satisfy \emph{the condition $HI_{r}$}. We will see that \emph{the condition $HI_r$} is a consequence of the $\cc_{r}$ condition (refer to Definition \ref{Crnotation}).

Let $R$ be an equidimensional, universally catenary, analytically unramified Noetherian local ring of dimension $d.$ Let $I=(a_1,\dots,a_g)$ be a proper $R$-ideal such that $\hgt(I) \geq 2.$ 

One of the principal tools we require to prove Theorem \ref{integralclosure} is that the integral closure of the powers of $I$ specializes when going modulo a general element $x$ (Notation \ref{notation}(f)). This was proved for $\ov{I}$ in \cite[Theorem 2.1]{hongUlrich} and the following theorem and its corollary generalizes it for higher powers $\ov{I^r},r\geq 2$.
\begin{Theorem} \label{specialprime}
	Let $\HH^2_J(\ov{\R'}(I))_{r-1} = 0$, for some $r$ and let $x = \sum_{i=1}^{g} z_ia_i \in R'.$ Then 
	\[ \ov{\left(\frac{(I')^r+(x)}{(x)}\right)} = \frac{\ov{(I')^r}+(x)}{(x)}. \]
\end{Theorem}

\begin{proof}
	Denote $R'/(x)$ by $S$ and $I'/(x)$ by $I'S.$ Let 
	\[ \R'(I'S)= \bigoplus_{n\in \Z} \left(\frac{(I')^n+(x)}{(x)}\right)t^n \]
	be the extended Rees ring of $I'S$ and $\ov{\R'}(I'S)$ be the integral closure of $\R'(I'S)$ in $S[t, t^{-1}].$ Let $J'=(t^{-1}, I't).$ Consider the natural map 
	\[ \phi \colon \ov{\R'}(I')/(xt)\ov{\R'}(I') \longrightarrow \ov{\R'}(I'S). \]
	Denote the kernel and cokernel of $\phi$ by $K$ and $C$ respectively. Our goal is to prove $C_r = 0.$ Following the arguments as in the proof of \cite[Theorem 2.1]{hongUlrich}, observe that $C_n \hookrightarrow \HH^2_{J'}(\ov{\R'}(I'))_{n-1}$, for all $n.$ As $R'$ is a flat $R$-extension, it follows that $\HH^2_{J'}(\ov{\R'}(I')) \simeq \HH^2_J(\ov{\R'}(I)) \otimes_R R'.$ Since, $\HH^2_J(\ov{\R'}(I))_{r-1} = 0$, we get $C_r=0$ and hence the result follows.
\end{proof}

\begin{Corollary} \label{special}
	Let $\HH^2_J(\ov{\R'}(I))_{r-1} = 0$, for some $r$ and let $x = \sum_{i=1}^{g}z_ia_i \in R''.$ Then
	\[ \ov{\left(\frac{(I'')^r+(x)}{(x)}\right)} = \frac{\ov{(I'')^r}+(x)}{(x)}. \]
\end{Corollary}
\begin{proof}
	As $R''$ is a flat $R$-extension, the proof follows using the same arguments as in the proof of Theorem \ref{specialprime}.
\end{proof}
We next prove  couple of lemmas required in the proof of Theorem \ref{integralclosure}. These lemmas describe the effect of going modulo the general element $x$ on the condition $\cc_{r}$.

\begin{Lemma} \label{pass}
	Suppose the normal filtration $\f = \{\ov{I^n}\}$ satisfies $\cc_{r}$ and let $S=R'/(x), I'S=I'/(x).$ Then the normal filtration $\f'=\{\ov{(I'S)^n} \}$ also satisfies $\cc_{r}$.	
\end{Lemma}

\begin{proof}
	Consider the natural map 
	\[ \phi \colon \ov{\R'}(I')/(xt)\ov{\R'}(I') \longrightarrow \ov{\R'}(I'S). \]
	Denote the kernel and cokernel of $\phi$ by $K$ and $C$ respectively. In \cite[Theorem 2.1]{hongUlrich}, the authors prove that $K = \HH^0_{J'}(K)$ and $C = \HH^0_{J'}(C)$, where $J'=(t^{-1},I't).$ This implies that $\HH^i_{J'}(K)=0$ and $\HH^i_{J'}(C)=0$, for all $i \ge 1.$ Consider the exact sequence 
	\[ 0 \longrightarrow K \longrightarrow \ov{\R'}(I')/(xt)\ov{\R'}(I') \overset{\phi}\longrightarrow \im(\phi) \longrightarrow 0. \]
	We obtain $\HH^i_{J'}(\ov{\R'}(I')/(xt)\ov{\R'}(I')) \simeq \HH^i_{J'}(\im(\phi))$ for all $i\ge 1.$ Then using the exact sequence 
	\[ 0 \longrightarrow \im(\phi) \longrightarrow \ov{\R'}(I'S) \longrightarrow C \longrightarrow 0, \]
	we obtain $\HH^i_{J'}(\ov{\R'}(I'S)) \simeq \HH^i_{J'}(\im(\phi)) \simeq \HH^i_{J'}(\ov{\R'}(I')/(xt)\ov{\R'}(I'))$ for all $i \ge 2.$ 
	
	Since $\hgt I \geq 1,$ we may assume that $a_1$ is a nonzerodivisor in $R$ and hence $x$ is a nonzerodivisor in $R'.$ Using the short exact sequence 
	\[ 0 \longrightarrow \ov{\R'}(I')(-1) \overset{xt}\longrightarrow \ov{\R'}(I') \longrightarrow \ov{\R'}(I')/(xt)\ov{\R'}(I') \longrightarrow 0, \]
	we get the following long sequence of local cohomology modules
	\[ \cdots \rightarrow \HH^i_{J'}(\ov{\R'}(I'))_j \rightarrow \HH^i_{J'}(\ov{\R'}(I')/(xt)\ov{\R'}(I'))_j \rightarrow \HH^{i+1}_{J'}(\ov{\R'}(I'))_{j-1} \rightarrow \cdots. \]
	As $R'$ is a flat $R$-module, it follows that $\HH^i_{J'}(\ov{\R'}(I')) \simeq \HH^i_J(\ov{\R'}(I)) \otimes_R R'.$ Thus $\HH^i_{J'}(\ov{\R'}(I'))_j=0$, for all $i,j$ such that $i+j=r+1$ and $0 \le i \le \hgt(I)$ (since $\f$ satisfies $\cc_r$). This implies that $\HH^i_{J'}(\ov{\R'}(I')/(xt)\ov{\R'}(I'))_j=0$ for all $i,j$ such that $i+j=r+1$ and $0 \le i \le \hgt(I)-1 = \hgt(I'S).$ Since $\HH^0_{J'}(\ov{\R'}(I'S))=0 = \HH^1_{J'}(\ov{\R'}(I'S))$ and $\HH^i_{J'}(\ov{\R'}(I'S)) \simeq \HH^i_{J'}(\ov{\R'}(I')/(xt)\ov{\R'}(I'))$ for all $2\le i \le 
\hgt(I'S)$, we are done. 
\end{proof}

Using similar arguments as above, one can prove the following corollary.
\begin{Corollary} \label{vanish}
	Suppose the normal filtration $\f=\{\ov{I^n}\}$ satisfies $\cc_{r}$. Then the normal filtration $\f''=\{\ov{(IT)^n} \}$ also satisfies $\cc_{r}$.
\end{Corollary}

\begin{Theorem} \label{integralclosure}
	Let $R$ be an equidimensional, universally catenary, and an analytically unramified Noetherian local ring of dimension $d.$ Let $I=(a_1,\dots,a_g)$ be an ideal generated by an $R$-regular sequence and $\f$ be the normal filtration of $I.$\\
	{\rm(1)} If $\hgt(I)=1,$ then {\em the condition $HI_r$} is true for all $r \ge 1.$ \\
	{\rm(2)} Let $\hgt(I) \ge 2$ and for some $r \ge 1$, suppose $\f$ satisfies $\cc_{r}$. Then $\f$ satisfies {\em the condition $HI_{r}$.}
\end{Theorem}

\begin{proof}
	Let $\hgt(I)=g$ and $I=(a_1,\dots,a_g)$, where $a_1,\dots a_g$ is an $R$-regular sequence. We apply induction on $g.$ Let $g=1$ and $I=(a).$ Write $\ov{R}$ for the integral closure of $R$ in its total ring of fractions. Since $a$ is a nonzerodivisor, we have $\ov{I} = (a)\ov{R} \cap R$ and $\ov{I^{n+r}} = (a^{n+r})\ov{R} \cap R.$ For all $n \ge 0$, we have
	\[ I^{n} \cap \ov{I^{n+r}} = (a^{n})R \cap (a^{n+r})\ov{R} = (a^{n})(R \cap (a^r)\ov{R}) = I^n\ov{I^r}. \]
	Thus the result holds when $g=1.$ Now define
	\begin{align*}
		G^{(r)}(I)=\frac{\R'(I)}{(t^{-r})} \simeq \bigoplus_{n\in \Z} \frac{I^n}{I^{n+r}} \ , && \ov{G^{(r)}}(I) = \frac{\ov{\R'}(I)}{(t^{-r})} \simeq\bigoplus_{n\in \Z} \frac{\ov{I^n}}{\ov{I^{n+r}}}.
	\end{align*}
	Consider the natural map 
	\begin{equation}\label{naturalses}
		0 \longrightarrow \ker(\phi^{(r)}) 
		\longrightarrow G^{(r)}(I) \otimes_R R/\ov{I^r} \simeq \bigoplus_{n\in \Z}\left(\frac{I^n}{I^n\ov{I^r}}\right) 
		\overset{\phi^{(r)}}\longrightarrow \bigoplus_{n\in\Z} \left(\frac{\ov{I^n}}{\ov{I^{n+r}}}\right) = \ov{G^{(r)}}(I).
	\end{equation}
	Then $K=\ker(\phi^{(r)}) = \bigoplus_{n\in\Z} (I^n \cap \ov{I^{n+r}})/I^n\ov{I^r}.$ Thus, to prove the theorem, it is enough to show $K=0.$
	
	Define $G^{(r)}(I'')=\bigoplus_{n\in\Z}(I'')^n/(I'')^{n+r}$ and $\ov{G^{(r)}}(I'')=\bigoplus_{n\in\Z}\ov{(I'')^n}/\ov{(I'')^{n+r}}.$ Observe that $R''$ is a faithfully flat extension of $R$, $\dim(R'') = \dim(R)$ and $I''$ is the ideal generated by an $R''$-regular sequence. It is easy to check that $R''$ is analytically unramified, equidimensional and universally catenary. Let $K''=K\otimes_R R''.$ Tensoring the equation (\ref{naturalses}) with $R''$, we get an exact sequence
	\begin{align}\label{defOfPhi''}
		0 \longrightarrow K'' \longrightarrow G^{(r)}(I'')\otimes_{R''} R''/\ov{(I'')^r} \overset{\phi''}\longrightarrow \ov{G^{(r)}}(I'').
	\end{align}
	
	Put $x=\sum_{i=1}^{g}a_iz_i.$ As $x$ is a nonzerodivisor in $R'',$ it follows that $\dim T=d-1$ and $\hgt(IT)=g-1.$ It is easy to check that the following  conditions are satisfied:
	(1) $T$ is analytically unramified, universally catenary and equidimensional (see, for example \cite[Theorem b]{hochsterGeneralGradeReductions}), and
	(2) $IT$ is an ideal generated by a $T$-regular sequence.
	
	Using Corollary \ref{vanish}, it follows that $\HH^i_{J''}(\ov{\R'}(IT))_j=0$ for all $i,j$ such that $i+j=r+1$ and $0 \le i \le \hgt(IT).$ This implies that $\HH^i_{J''/(x)}(\ov{\R'}(IT))_j=0$ for all $i,j$ such that $i+j=r+1$ and $0 \le i \le \hgt(IT).$ Therefore, we can now apply induction hypothesis to $T$ and $IT.$ Let $G^{(r)}(IT) = \bigoplus_{n\in \Z} (IT)^n/(IT)^{n+r}$ and $\ov{G^{(r)}}(IT) = \bigoplus_{n\in\Z} \ov{(IT)^n}/\ov{(IT)^{n+r}}.$ By induction hypothesis, we get the exact sequence
	\[ 0 \rightarrow  G^{(r)}(IT) \otimes_{T} T/\ov{(IT)^r} \overset{\widetilde{\phi}}\longrightarrow  \ov{G^{(r)}}(IT). \]
	
	Observe that $t^{-1},xt$ is a regular sequence in $\R'(I'').$ This implies that $t^{-r},xt$ is a regular sequence in $\R'(I'')$ and hence image of $x$ in $G^{(r)}(I'')$ is a regular element in $G^{(r)}(I'').$ Similarly, the image of $x$ in $\ov{G^{(r)}}(I'')$ is also regular element in $\ov{G^{(r)}}(I'').$ Consider 
	\[ \left(G^{(r)}(IT) \right)_n=\frac{(IT)^n}{(IT)^{n+r}} = \frac{(I'')^n+(x)}{(I'')^{n+r}+(x)} \simeq \frac{(I'')^n}{(I'')^{n+r}+((x)\cap (I'')^n)}. \]
	Now, 
	\[ \left( \frac{G^{(r)}(I'')}{((x)+(I'')^{r+1})} \right)_n \simeq \frac{(I'')^n}{((x)(I'')^{n-1}+(I'')^{n+r})}. \]
	As image of $x$ is a nonzerodivisor in $G^{(r)}(I'')$, $(x)(I'')^{n-1} = (x)\cap (I'')^n$, for all $n.$ This implies that $G^{(r)}(IT) \simeq G^{(r)}(I'')/((x)+(I'')^{r+1}).$ Also, $\ov{\R'}(I'')$ is a finite $\R'(I'')$-module implies that $\ov{G^{(r)}}(I'')$ is a finite $G^{(r)}(I'')$-module.

	Consider the short exact sequence of finite graded modules over $G^{(r)}(I'')$
	\[ 0 \rightarrow K'' \rightarrow G^{(r)}(I'')\otimes_{R''}R''/\ov{(I'')^r} \overset{\phi''}\longrightarrow \im(\phi'') \rightarrow 0. \]
	As $(x)+(I'')^{r+1}$ is regular on $\im(\phi'')$, on tensoring the sequence with $G^{(r)}(IT)$ over $G^{(r)}(I'')$, we get the short exact sequence (using \cite[Proposition 1.1.4]{brunsHerzog})
	\[ 0 \rightarrow G^{(r)}(IT)\otimes_{G^{(r)}(I'')}K'' \rightarrow G^{(r)}(IT)\otimes_{R''}R''/\ov{(I'')^r}\rightarrow G^{(r)}(IT)\otimes_{G^{(r)}(I'')}\im(\phi'') \rightarrow 0. \]
	Notice that $ G^{(r)}(IT) \otimes_{T} T/\ov{(IT)^r}\simeq   G^{(r)}(IT) \otimes_{T}\left(R''/(x)\otimes_{R''} R''/\ov{(I'')^r}\right)\simeq G^{(r)}(IT)\otimes_{R''}R''/\ov{(I'')^r}$ (where the first isomorphism is due to Corollary \ref{special}). This gives a commutative diagram with exact rows
	\[\xymatrixrowsep{6mm} \xymatrixcolsep{6mm}
	\xymatrix{
		0 \ar[r] &
		G^{(r)}(IT) \otimes_{G^{(r)}(I'')} K'' \ar[r] & 
		G^{(r)}(IT) \otimes_{R''}R''/\ov{(I'')^r} \ar[r]^{\phi''} \ar[d]_\alpha^\simeq & 
		G^{(r)}(IT) \otimes_{G^{(r)}(I'')} \im(\phi'') \ar[r] \ar[d]^\psi &
		0\\
		& 0 \ar[r] &  G^{(r)}(IT) \otimes_{T} T/\ov{(IT)^r} \ar[r]^{\widetilde{\phi}}  &  \ov{G^{(r)}}(IT) 
	} \]
For ease of notation, in the above diagram, we refer to $\phi''$ as the map induced from $\phi''$ of \eqref{defOfPhi''}. We define $\psi$ as follows. 
	%Observe that $\im(\phi'')=\bigoplus_{n\in\Z}(I'')^n+\ov{(I'')^{n+r}}/\ov{(I'')^{n+r}}.$
	As $\im(\phi'')$ is a $G^{(r)}(I'')$-module and $G^{(r)}(IT) \simeq G^{(r)}(I'')/(x+(I'')^{r+1})$, we get 
	\[ G^{(r)}(IT) \otimes_{G^{(r)}(I'')} \im(\phi'') \simeq \im(\phi'')/(x+(I'')^{r+1})\im(\phi''). \]
	
	Define $\delta: \im(\phi'') \rightarrow \ov{G^{(r)}}(IT)$ by $\delta(a) = a+\ov{(x+(I'')^{r+1})}.$ Then $\delta$ is a well defined map and as $(x+(I'')^{r+1})\im(\phi'') \subseteq \ker(\delta)$, there exists a map $\psi: G^{(r)}(IT) \otimes_{G^{(r)}(I'')}\im(\phi'') \rightarrow \ov{G^{(r)}}(IT).$ Since this map is induced by the natural maps, the above diagram is commutative. We now show that $G^{(r)}(IT) \otimes_{G^{(r)}(I'')}K''=\ker(\phi'') =0.$ Let $a \in \ker(\phi'').$ Then $\phi''(a) = 0$ implies that $(\psi \circ \phi'')(a)=0 = (\widetilde{\phi} \circ \alpha)(a).$ As $\widetilde{\phi}$ is an injective map, we get $\alpha(a)=0$ and hence $a=0.$ Using graded Nakayama Lemma, we get $K=0.$
\end{proof}

Using the above theorem, we recover Huneke-Itoh intersection theorem (\cite[Theorem 1]{itoh88}) for ideals with height 2.

\begin{Corollary}
	Let $R$ be an equidimensional, universally catenary, and analytically unramified Noetherian local ring. Let $I$ be an ideal generated by an $R$-regular sequence such that $\hgt I =2.$ Then for all $n \geq 0,$ 
	\[ I^n \cap \overline{I^{n+1}} = I^n \overline{I}. \]
\end{Corollary}

\begin{proof} 
	We need to show that the filtration $\f=\{\ov{I^n} \}$ satisfies {\em the condition $HI_1.$} Using \cite[Proposition 13]{itoh92}, it follows that $\HH^0_J(\ov{\R'}(I))=0=\HH^1_J(\ov{\R'}(I)).$ Since $\HH^2_J(\ov{\R'}(I))_0=0$ from \cite[Theorem 2]{itoh88}, using Theorem \ref{integralclosure}, we are done. 
\end{proof}
	
If $I$ is generated by a regular sequence in  Noetherian local ring, then \emph{the condition $HI_1$} on $I$ is the Huneke-Itoh Intersection Theorem. But \emph{the condition $HI_2$} may not be satisfied by the normal filtration of $I.$ We illustrate this in the next example.

\begin{Example}\label{exampleHunekeHuckaba}{\rm 
	Let $R$ be a $3$-dimensional regular local ring, $I$  an $R$-ideal generated by an $R$-regular sequence of length 3 and suppose $I^{n}\cap \overline{I^{n+2}}=\overline{I^2}I^{n}$, for all $n.$ Then by the Brian\c{c}on-Skoda Theorem, $\overline{I^{w+3}}\subseteq I^{w+1}$ for $w\geq 0.$ Thus for $n\geq 1$, we have $\overline{I^{n+2}}\subseteq I^{n}.$ Therefore, $\overline{I^{n+2}}=I^{n}\cap \overline{I^{n+2}}=\overline{I^2}I^{n}.$ This shows that the normal reduction number is at most $2.$  However, Huckaba and Huneke \cite[Theorem 3.11]{HunekeHuckaba99}, constructed a system of parameters  in a 3-dimensional regular local ring which has normal reduction number greater than 2. We describe this example.
		
	Let $k$ be a field of characteristic not equal to 3. Consider the ring $R=k[[x,y,z]]$ and an ideal $I=(x^4, z(y^3+z^3), y(y^3+z^3)+z^5)$ generated by an $R$-regular sequence. Then $\ov{I}=N+\m^5$, where $N=(x^4,x(y^3+z^3),y(y^3+z^3),z(y^3+z^3))$ and $\m=(x,y,z).$ It is proved in \cite[Theorem 3.11]{HunekeHuckaba99} that $\ov{I}$ is a height 3 normal ideal and $\ov{G}(I)$ is not Cohen-Macaulay. Since $\ov{G}(I)$ is not Cohen-Macaulay, using \cite[Proposition 3]{itoh88}, normal reduction number of $I$ is greater than 2. Thus, $\f=\{\ov{I^n}\}$ does not satisfy \emph{the condition $HI_2.$}
	
	Interestingly, in \cite[Theorem 3.11]{HunekeHuckaba99}, the authors also prove that $\HH^3_J(\ov{\R'}(I))_0 \neq 0$. Thus  we believe that if the cohomological  conditions in Theorem \ref{integralclosure} fail, then it is less likely that the filtration $\f=\{ \overline{I^n}\}$ will satisfy {\em the condition $HI_2$}.
}\end{Example}

\begin{Remarks}{\rm 
		(1) In (\cite[Lemma 14]{itoh92}, \cite[Theorem 2]{ooishi}), the author proves that if $\HH^i_{\R(I)_+}(\ov{\R'}(I))_j=0$ for all $i,j$ such that $i+j=r+1$, then $\ov{r}(I) \leq r.$ This implies that \emph{the condition $HI_{r}$} is true as for all $n \ge 0$, $\ov{I^{n+r}}\cap I^{n} = I^n\ov{I^r}\cap I^{n} = I^{n}\ov{I^r}.$

\noindent (2) The exact sequence (\ref{itohseq}) implies that if $\HH^d_{\R(I)_+}(\ov{\R'}(I))_{r-d+1}=0$,  then $\HH^d_J(\ov{\R'}(I))_{r-d+1}=0$ but the converse may not be true. Thus if $\hgt (I)=d,$ then the assumptions made in Theorem \ref{integralclosure} are weaker than assuming that $\mathcal{F}$ satisfies $\hc_{r}$. %$\HH^i_{\R(I)_+}(\ov{\R'}(I))_{j} = 0$ for all $i,j$ such that $i+j=r+1$ and $0 \le i \le \hgt(I).$ 
We illustrate this with an example where $\HH^d_J(\ov{\R'}(I))_{j}=0$ but $\HH^d_{\R(I)_+}(\ov{\R'}(I))_j \not= 0.$
}\end{Remarks}

\begin{Example}{\rm
	Let $(R,\m)$ be a 2-dimensional Cohen-Macaulay, analytically unramified local ring and $I$ be a complete intersection such that $\ov{r}(I)=2.$ By \cite[Proposition 3]{itoh88} and \cite[Corollary 3.8]{marleyThesis}, we get $\ov{n}(I)=0.$ Hence $\ov{P}_I(0)-\ov{H}_I(0) = \ell(\HH^2_{\R(I)_+}(\ov{\R'}(I))_0) \not=0$ using the Difference Formula. But $\HH^2_{J}(\ov{\R'}(I))_0 = 0$ from \cite[Theorem 2]{itoh88}.
	
	Let $R=\mathbb{C}[[X,Y,Z]]/(X^3+Y^3+Z^3).$ Then $R$ is a 2-dimensional Cohen-Macaulay, analytically unramified local ring. Let $x,y,z$ denote the images of $X,Y$ and $Z$ in $R$ respectively and let $I=(y,z).$ Then $I$ is a complete intersection such that $I$ is a reduction of the maximal ideal $\m.$ As $G(\m)=\mathbb{C}[X,Y,Z]/(X^3+Y^3+Z^3)$ is reduced, it follows that $\overline{I^n} = \overline{\m^n} = \m^n$ for all $n.$ We now claim that $\ov{r}(I) = r(\m) = 2.$ If $\m I = \m^2$, then Abhaynkar's equality is true and so we must have $v(\m)-\dim R +1 = e(\m)$, where $v(\m)$ denotes the embedding dimension. But $v(\m)-\dim R+1 = 3-2+1 = 2 \neq e(\m).$ Hence $\m I \neq \m^2.$ As $\m^2 I = \m^3$, the claim holds. Since $G(\m)=\ov{G}(I)$ is Cohen-Macaulay, we get $\ov{n}(I)=0$ and using the Difference Formula, it follows that $\ov{P}_I(0)-\ov{H}_I(0) = \ell(\HH^2_{\R(I)_+}(\ov{\R'}(I))_0) \not=0.$ But from Itoh's theorem, $\HH^2_{J}(\ov{\R'}(I))_0=0.$
}\end{Example}

\section{Local Cohomology modules and normal reduction number}

In view of results of Itoh \cite{itoh92}, we know that there is a strong relationship between the  the normal reduction number and the normal Hilbert coefficients (particularly  $\ov{e}_2(I)$ and $\ov{e}_1(I)$). Itoh shows that $\ov{e}_2(I)=\ov{e}_1(I)-\ell(\ov{I}/I)$ is equivalent to $\ov{r}(I)\leq 2$. The main purpose of this section is to generalize this idea and check if the generalized formulas for $\ov{e}_2(I)$ leads to bounds on the normal reduction number. This section also builds the necessary tools required in the next section.

The setting of this section is same as that of Section \ref{normafiltrationHI}. The following lemma is a generalization of \cite[Lemma 14, Corollary 16]{itoh92}.

\begin{Lemma}\label{reduNumlocalCohom}
	Let $(R,\m)$ be a $d$-dimensional analytically unramified, Cohen Macaulay local ring with $d\geq 2$ and let $I$ be a parameter ideal. For some $k \geq 2$, suppose the filtration $\f=\{\ov{I^n}\}$ satisfies \emph{the condition $HI_p$} for all $p \leq k-2$. Then $\ov{r}(I)\leq k-1$ if and only if $\f$ satisfies $\hc_{k-1}$.
\end{Lemma} 

\begin{proof}
	Suppose $\HH^i_{\R(I)_+}(\ov{\R'}(I))_j=0$ for all $i,j$ such that $i+j = k$ and $2 \leq i \leq d$. Then using \cite[Lemma 14]{itoh92}, it follows that $\ov{r}(I)\leq k-1$.
	
	Conversely, suppose that $\ov{r}(I)\leq k-1$. As the ideal $I$ satisfies \emph{the condition $HI_p$} for all $p \leq k-2$, it follows that the associated graded ring $\ov{G}(I)$ and hence the extended Rees algebra $\ov{\mathcal{R}'}(I)$ is Cohen-Macaulay. Thus $\HH^i_{\R(I)_+}(\ov{\R'}(I))=0$ for $0\leq i\leq d-1$. Also, $\HH^d_{\R(I)_+}(\ov{\R'}(I))_{k-d}=0$ (by \cite[Lemma 15]{itoh92}) and hence the result.
\end{proof}

Note that in the above lemma, if $\f$ satisfies $\hc_{k-1}$, then $\f$ satisfies \emph{the condition $HI_{k-1}$} (using Theorem \ref{integralclosure} and \cite[Proposition 13]{itoh92}). Thus if $\ov{r}(I)\leq k-1$, then \emph{the condition $HI_{k-1}$} is satisfied, which is an expected conclusion.

The following lemma is a generalization of \cite[Proposition 17]{itoh92}.

\begin{Lemma}\label{reductionNumberpassGenExt}
	Let $(R,\m)$ be a $d$-dimensional analytically unramified, Cohen Macaulay local ring with $d\geq 2$ and $I$ be a parameter ideal. For some $k \geq 2$, suppose the normal filtration $\f$ satisfies $\cc_{r}$ for $2 \leq r \leq k-2$. If $\ov{r}(IT) \leq k-1$, then $\ov{r}(I) \leq k-1.$ 
\end{Lemma}

\begin{proof}
	Corollary \ref{vanish}, Theorem \ref{integralclosure} and \cite[Theorem 1]{itoh88} shows that $IT$  satisfies \emph{the condition $HI_p$} for all $p \leq k-2$, in the ring $T$. Since $\ov{r}(IT)\leq k-1$, using Lemma \ref{reduNumlocalCohom} and \cite[Lemma 14]{itoh92} it follows that
	\begin{align*}
		\HH^i_{\R(IT)_+}(\ov{\R '}(IT))_j=0 \text{ for } 2\leq i\leq d-1 \text{ and } i+j \geq k.
	\end{align*}
	
	Consider the exact sequence 
	\begin{align*}
		0\rightarrow K\rightarrow \frac{\ov{\R'}(I'')}{(xt)\ov{\R'}(I'')}\rightarrow \ov{\R'}(IT)\rightarrow C\rightarrow 0.
	\end{align*}
	Using the ideas as in the the proof of \cite[Theorem 2.1]{hongUlrich}, we have $K=\HH^0_{\R(I'')_+}(K), C=\HH^0_{\R(I'')_+}(C)$ and $\HH^i_{\R(I'')_+}(K)=\HH^i_{\R(I'')_+}(C)=0$ for $i\geq 1$ (Notice that the support of the local cohomology modules in here and in \cite[Theorem 2.1]{hongUlrich} are different, but the proof is essentially the same). Thus as in the proof Lemma \ref{pass}, we have
	\begin{align}\label{sec5eq1}
		\HH^i_{\R(I'')_+}\left(\frac{\ov{\R'}(I'')}{(xt)\ov{\R'}(I'')} \right)_j
		\cong\HH^i_{\R(IT)_+}(\ov{\R '}(IT))_j=0 \text{ for }2\leq i\leq d-1 \text{ and } i+j \geq k.
	\end{align}
	The set of isomorphisms also require the change of rings formula in local cohomology. Now consider the short exact sequence
	\begin{align*}
		0\rightarrow (\ov{\R'}(I''))(-1)\overset{xt}\longrightarrow \ov{\R'}(IR'')\rightarrow \left(\frac{\ov{\R'}(I'')}{(xt)\ov{\R'}(I'')} \right)\rightarrow 0
	\end{align*}
	from which we have a long exact sequence of local cohomology modules
	\begin{multline*}
		\cdots\rightarrow 
		\HH^i_{\R(I'')_+}(\ov{\R'}(I''))_{j-1}\rightarrow \HH^i_{\R(I'')_+}(\ov{\R'}(I''))_j\rightarrow \HH^i_{\R(I'')_+}\left(\frac{\ov{\R'}(I'')}{(xt)\ov{\R'}(I'')} \right)_j\rightarrow\cdots.
	\end{multline*}
	Since $R\rightarrow R''$ is a faithfully flat extension, using the hypothesis and \cite[Proposition 13]{itoh92} it follows that 
	\[ \HH^i_{\R(I'')_+}(\ov{\R'}(I''))_j \cong  \HH^i_{J''}(\ov{\R'}(I''))_j=0 \] 
	for $3 \leq i+j=r+1\leq k-1$ and $2\leq i\leq d-1$.  For $0\leq i\leq d-1$ and $i+j=k$, the long exact sequence of local cohomology modules gives
	\begin{align*}
		0=\HH^i_{\R(I'')_+}(\ov{\R'}(I''))_{j-1}=\HH^i_{\R(I'')_+}(\ov{\R'}(I''))_j
	\end{align*}
	where the first equality follow from the fact that $i+j-1=k-1$ and the second equality follows from \eqref{sec5eq1}. For $i=d$ and $j\geq k-d$, we have from the long exact sequence
	\begin{multline*}
		\cdots\rightarrow 0=\HH^{d-1}_{\R(I'')_+}\left(\frac{\ov{\R'}(I'')}{(xt)\ov{\R'}(I'')} \right)_{j+1} \rightarrow 
		\HH^d_{\R(I'')_+}(\ov{\R'}(I''))_{j}
		\overset{xt}\longrightarrow \HH^d_{\R(I'')_+}(\ov{\R'}(I''))_{j+1} \rightarrow\cdots
	\end{multline*}
	This implies that the map $\HH^d_{\R(I'')_+}(\ov{\R'}(I''))_j \overset{xt} \longrightarrow \HH^d_{\R(I'')_+}(\ov{\R'}(I''))_{j+1}$ is injective for all $j \geq k-d.$ As some power of $xt$ annihilates these local cohomology modules, we have $\HH^d_{\R(I'')_+}(\ov{\R'}(I''))_{k-d}=0$. Thus, $\HH^i_{\R(I)_+}(\ov{\R'}(I))_j=0$ for all $i,j$ such that $i+j = k$ and $2 \leq i \leq d$ and hence $\ov{r}(I)\leq k-1$ by Lemma \ref{reduNumlocalCohom}.
\end{proof}

\begin{Proposition}
	Let $R$ be a $2$-dimensional  analytically unramified, Cohen Macaulay local ring and $I$ be a parameter ideal. For some $k \geq 2,$ if  
	\[ \ov{e}_2(I)=(k-2)\ov{e}_1(I)-\sum_{i=0}^{k-3}(k-2-i) \ \ell\left(\frac{\ov{I^{i+1}}}{I\ov{I^i}}\right),  \]
	then $\ov{r}(I)\leq k-1.$
\end{Proposition}

\begin{proof}
	By \cite[Remark 4.2, 4.3]{huneke1987} we have, for all large $n$,
	\begin{align*}
		\ell\left(\frac{R}{\ov{I^{n+1}}}\right)
		=\ell\left(\frac{R}{I}\right){n+2\choose 2}-\ov{e}_1(I) {n+1\choose 1} + \ov{e}_2(I),\\
		\ov{e}_1(I)=\sum_i\ell\left(\frac{\ov{I^{i+1}}}{I\ov{I^i}}\right) 
		\ \text{ and } \
		\ov{e}_2(I)=\sum_i i \cdot \ell\left(\frac{\ov{I^{i+1}}}{I\ov{I^i}}\right).	
	\end{align*}
	Thus if  $\ov{e}_2(I)=(k-2)\ov{e}_1(I)-\sum_{i=0}^{k-3}(k-2-i)\ell\left(\dfrac{\ov{I^{i+1}}}{I\ov{I^i}}\right)$, then by the above equalities we have $\sum_{i\geq k-1} (i-k+2)  \ell\left(\dfrac{\ov{I^{i+1}}}{I\ov{I^i}}\right)=0$ which implies that $\ov{r}(I)\leq k-1$.
\end{proof}

The following theorem generalizes the result of Itoh mentioned above.

\begin{Theorem}\label{e2theorem}
	Let $(R,\m)$ be a $d$-dimensional analytically unramified, Cohen Macaulay local ring with $d\geq 2$ and let $I$ be a parameter ideal. For some $k \geq 2,$ suppose that $\f$ satisfies $\cc_{r}$ for $2 \leq r \leq k-2$. If $\ov{e}_2(I)=(k-2)\ov{e}_1(I)-\sum_{i=0}^{k-3}(k-2-i)\ell\left(\dfrac{\ov{I^{i+1}}}{I\ov{I^i}}\right)$, then $\ov{r}(I)\leq k-1$.
\end{Theorem}

\begin{proof}
	We induct on the dimension $d$. If $d=2$, then the result is true by the previous proposition. Now suppose $d>2$. Pass to the ring $T.$ Since $\ov{e}_1(I)=\ov{e}_1(IT)$ and $\ov{e}_2(I)=\ov{e}_2(IT)$ (\cite[Corollary 8]{itoh92}), due to Theorem \ref{specialprime}, the formula for $\ov{e}_2(IT)$ is preserved. Now by Corollary \ref{vanish} and induction hypothesis, $\ov{r}(IT)\leq k-1$. Now $\ov{r}(I)\leq k-1$ follows from Lemma \ref{reductionNumberpassGenExt}.
\end{proof}

\section{Vanishing criterion for normal Hilbert coefficients}

Itoh, in \cite{itoh92}, proposed the following conjecture relating the vanishing of $\overline{e}_3(I)$ with the normal reduction number of the ideal $I$.
\begin{Conjecture}
	Let $(R,\m)$ be an analytically unramified Gorenstein local ring of dimension $d \geq 3.$ Let $I$ be a parameter ideal. Then $\overline{e}_3(I) = 0$ if and only
	if $\overline{r}(I) \leq 2.$
\end{Conjecture}
Itoh solved the conjecture when $\ov{I} = \m.$ Using  \emph{the condition $HI_r$}, we prove an analogue of Itoh's theorem for vanishing of $\ov{e}_k(I)$, for $k\leq d=\dim R$ and $d \geq 3.$ We first prove some preliminary results.

\begin{Proposition}\label{HSpoly}
	Let $R$ be a $d$-dimensional analytically unramified, Cohen-Macaulay local ring and $I$ be a parameter ideal. Let $\f=\{ \ov{I^n} \}$ satisfy \emph{the condition $HI_p$} for all $p \leq k-2$ and let $d \geq k-1.$ 
	\begin{enumerate}[$(a)$]
		\item 	Then for all $n \geq k-2,$
		\begin{equation} \label{poly}
			\ell\left(\frac{R}{\ov{I^{n+1}}}\right)
			\leq \ell\left(\frac{R}{I}\right)\binom{n+d}{d} - \alpha_1(\f)\binom{n+d-1}{d-1} + \cdots + (-1)^{k-1} \alpha_{k-1}(\f)\binom{n+d-(k-1)}{d-(k-1)}
		\end{equation}
		where $\displaystyle \alpha_j(\f) = \sum_{i=j-1}^{k-2} \binom{i}{j-1}\ell(\ov{I^{i+1}}/I\ov{I^i})$, for all $j=1,\ldots, k-1.$ The equality holds in the equation {\rm(\ref{poly})} if and only if $\ov{r}(I) \leq k-1.$ In this case, $\ov{G}(I)$ is Cohen-Macaulay.
		\smallskip
		\item The equality $\ov{e}_1(I) = \sum_{j=0}^{k-2} \ell(\ov{I^{j+1}}/I\ov{I^j})$ holds if and only if $\ov{r}(I)\leq k-1$. In this case, $\ov{G}(I)$ is Cohen-Macaulay.
	\end{enumerate}
\end{Proposition}	

\begin{proof}
	$(a)$: Note that for all $n \geq k-2,$
	\begin{align*}
		\ell\left(\frac{R}{\ov{I^{n+1}}}\right)
		&= \ell\left(\frac{R}{I^n}\right) + \ell\left(\frac{I^n}{I^n\ov{I}}\right) - \ell\left(\frac{\ov{I^{n+1}}}{I^n\ov{I}}\right)\\
		&= \ell\left(\frac{R}{I}\right){n+d-1\choose d} + \ell\left(\frac{I^n}{I^n\ov{I}}\right) - \sum_{i=1}^{k-2} \ell\left(\frac{I^{n-i}\ov{I^{i+1}}}{I^{n-i+1}\ov{I^{i}}}\right) - \ell\left(\frac{\ov{I^{n+1}}}{I^{n-(k-2)}\ov{I^{k-1}}}\right).
	\end{align*}
	Since $I^n/I^{n+1}\otimes_R R/\ov{I} \simeq I^n/I^n\ov{I}$ and $I^n/I^{n+1}$ is a free $R/I$-module, it follows that 
	{\small \begin{align}\label{HSpolyeq1}
		\ell\left(\frac{R}{\ov{I^{n+1}}}\right)
		&= \ell\left(\frac{R}{I}\right){n+d-1\choose d} + {n+d-1\choose d-1}\ell\left(\frac{R}{\ov{I}}\right) - \sum_{i=1}^{k-2} \ell\left(\frac{I^{n-i}\ov{I^{i+1}}}{I^{n-i+1}\ov{I^{i}}}\right) - \ell\left(\frac{\ov{I^{n+1}}}{I^{n-(k-2)}\ov{I^{k-1}}}\right).
	\end{align}}
	As $\f$ satisfies \emph{the condition $HI_p$} for all $p \leq k-2$, from \cite[Lemma 20]{GMV} it follows that for all $i=1,\ldots,k-2$ and $n \ge i$,
	\[ \frac{I^{n-i}\ov{I^{i+1}}}{I^{n-i+1}\ov{I^{i}}}\simeq \frac{I^{n-i}}{I^{n-i+1}}\bigotimes_R\frac{\ov{I^{i+1}}}{I\ov{I^i}}. \]
	Moreover, $I^{n-i}/I^{n-i+1}$ is a free $R/I$-module and so 
	\[ \ell\left(\frac{I^{n-i}\ov{I^{i+1}}}{I^{n-i+1}\ov{I^i}}\right) 
	= \ell\left(\frac{\ov{I^{i+1}}}{I\ov{I^i}}\right){n+d-i-1\choose d-1} 
	= \ell\left(\frac{\ov{I^{i+1}}}{I\ov{I^i}}\right) \left[\sum_{j=1}^{i+1} (-1)^{j-1} \binom{i}{j-1}\binom{n+d-j}{d-j} \right] \]
	where the last equality follows from \cite[Lemma 22]{GMV}.
	Hence for $n \ge k-2$,
	\begin{align}
		\ell\left(\frac{R}{\ov{I^{n+1}}}\right)
		&= \ell\left(\frac{R}{I}\right)\binom{n+d}{d} - \ell\left(\frac{\ov{I}}{I}\right)  \binom{n+d-1}{d-1} \nonumber\\
		&- \sum_{i=1}^{k-2} \ell\left(\frac{\ov{I^{i+1}}}{I\ov{I^i}}\right) \left[\sum_{j=1}^{i+1} (-1)^{j-1} \binom{i}{j-1}\binom{n+d-j}{d-j} \right] - \ell\left(\frac{\ov{I^{n+1}}}{I^{n-(k-2)}\ov{I^{k-1}}}\right).\label{HSpolyeq2} \\
		&= \ell\left(\frac{R}{I}\right)\binom{n+d}{d} - \alpha_1(\f)\binom{n+d-1}{d-1} + \cdots + (-1)^{k-1} \alpha_{k-1}(\f)\binom{n+d-(k-1)}{d-(k-1)} \nonumber\\
		&- \ell\left(\frac{\ov{I^{n+1}}}{I^{n-(k-2)}\ov{I^{k-1}}}\right)\nonumber
	\end{align}
	where $\displaystyle \alpha_j(\f) = \sum_{i=j-1}^{k-2} \binom{i}{j-1} \ell(\ov{I^{i+1}}/I\ov{I^i})$, for all $j=1,\ldots,k-1.$
	Hence equation (\ref{poly}) holds. 
	
	If equality holds in the equation (\ref{poly}), then  $\ov{I^{n+1}}=I^{n-(k-2)}\ov{I^{k-1}}$ for all $n \geq k-2$ which implies that $\ov{r}(I) \leq k-1.$ Observe that as $\mathcal{F}$ satisfies \emph{the condition $HI_p$} for all $p \leq k-2$ and as $\ov{r}(I) \leq k-1$, using Valabrega-Valla criterion (or \cite[Theorem 25]{GMV}), it follows that $\ov{G}(I)$ is Cohen-Macaulay. Conversely, if $\ov{r}(I) \leq k-1$, then $\ov{I^{n+1}}=I^{n-(k-2)}\ov{I^{k-1}}$ for all $n \geq k-2$ and hence equality holds in the equation (\ref{poly}).
	\smallskip
	
	$(b)$: 	Notice that 
	\begin{align*}
		\sum_{j=0}^{k-2}\ell\left(\frac{\ov{I^{j+1}}}{I\ov{I^j}}\right)
		=\sum_{j=0}^{k-2}\ell\left(\frac{\ov{I^{j+1}}}{I\cap \ov{I^{j+1}}} \right)
		=\sum_{j=0}^{k-2}\ell\left(\frac{\ov{I^{j+1}}+I}{I} \right)
	\end{align*}
	where the first equality is due to the fact that $\f$ satisfies \emph{the condition $HI_p$} for all $p \leq k-2.$ Now \cite[Corollary 4.8]{huckabaMarley} and \cite[Theorem 25]{GMV} show that the equality $\ov{e}_1(I)=\sum_{j=0}^{k-2}\ell(\ov{I^{j+1}}+I/I)$ holds if and only if $\ov{r}(I)\leq k-1$. Using same arguments as in part $(a)$, we get $\overline{G}(I)$ is Cohen-Macaulay.
	
	The converse statement can also be recovered from part $(a)$. If $\ov{r}(I)\leq k-1$, then equality holds in the equation \eqref{poly} and hence it follows that $\ov{e}_1(I)=\sum_{j=0}^{k-2}\ell(\ov{I^{j+1}}/I\ov{I^j})$.
\end{proof}

\begin{Theorem} \label{vanishingOfe_k}
	Let $(R,\m)$ be a $d$-dimensional analytically unramified, Cohen-Macaulay local 
	ring and let $I$ be a parameter ideal. Suppose $k \leq d$ and $d\geq 2.$  If $\HH^i_{\R(I)_+}(\ov{\R'}(I))_j=0$ for $3\leq i+j\leq k-1$ and $3\leq i\leq d-1$, then $\ov{e}_k(I)\geq 0$ and if $\ov{e}_k(I)=0$ then $\ov{I^{n+k-1}}\subseteq I^n$, for every $n\geq 0$.
\end{Theorem}

\begin{proof}
	The result is known for $k\leq 3$ (see\cite{itoh92}). But we prove it for the sake of completeness. We prove by induction on $d.$ Let $d=2.$ If $k=2$, then using the Difference Formula, 
	\[ \ov{e}_2(I) = \ov{P}_{I}(0) - \ov{H}_{I}(0) = \sum_{i=0}^2 (-1)^i \ell_R(\HH^i_{\R(I)_+}(\ov{\R'}(I))_0).\]
	From \cite[Proposition 13]{itoh92}, it follows that $\HH^0_{\R(I)_+}(\ov{\R'}(I)) = 0 = \HH^1_{\R(I)_+}(\ov{\R'}(I))$ and hence $\ov{e}_2(I) = \ell_R(\HH^2_{\R(I)_+}(\ov{\R'}(I))_0) \geq 0.$ If $\ov{e}_2(I)=0,$ then $\HH^2_{\R(I)_+}(\ov{\R'}(I))_0=0$ which implies that $\ov{I^{n+1}}\subseteq I^n$ for every $n\geq 0$, using \cite[Lemma 18]{itoh92}. Now let $d=2$ and $k=1.$ Using \cite[Theorem 4.7]{huckabaMarley}, 
	\[ \ov{e}_1 \geq \sum_{n \geq 1}\ell\left( \frac{\ov{I^n}}{I \cap \ov{I^n}} \right) \geq 0 \]
	and if $\ov{e}_1=0,$ then $\ov{I^n} = I \cap \ov{I^n}$ for all $n \geq 1.$ This implies that $\ov{I^n} \subseteq I^n$ for all $n.$ 
	
	Suppose $d \geq 3.$ If $k=d$, then using the Difference Formula, 
	\[ (-1)^d \ov{e}_k(I) = \ov{P}_{I}(0) - \ov{H}_{I}(0) = \sum_{i=0}^d (-1)^i \ell_R(\HH^i_{\R(I)_+}(\ov{\R'}(I))_0). \]
	Note that $\HH^0_{\R(I)_+}(\ov{\R'}(I)) = 0 = \HH^1_{\R(I)_+}(\ov{\R'}(I))$ and $\HH^2_J(\ov{\R'}(I))_0 \simeq \HH^2_{\R(I)_+}(\ov{\R'}(I))_0=0$ using \cite[Proposition 13]{itoh92} and \cite[Theorem 2]{itoh88}. By assumption and \cite[Proposition 13]{itoh92}, 
	\[ \HH^3_{\R(I)_+}(\ov{\R'}(I))_0 = \cdots = \HH^{d-1}_{\R(I)_+}(\ov{\R'}(I))_0=0 \]
	and hence $\ov e_d(I)= \ell(\HH^d_{\R(I)_+}(\ov{\R'}(I))_0) \geq 0$. If $\ov e_d(I)=0$, then $\HH^d_{\R(I)_+}(\ov{\R'}(I))_0=0$. This implies that $\ov{I^{n+d-1}}\subseteq I^n$ for every $n\geq 0$, using \cite[Lemma 18]{itoh92}. Now let $k < d.$  Using Corollary \ref{vanish}, $\HH^i_J(\ov{\R'}(I))_j=0$ for all $i,j$ such that $3 \leq i+j \leq k-1$ and $3 \leq i \leq d-1$ implies that $\HH^i_{J''}(\ov{\R'}(IT))_j=0$ for all $i,j$ such that $3 \leq i+j \leq k-1$ and $3 \leq i \leq d-2.$ Since $T$ is an analytically unramified, Cohen-Macaulay local ring and $IT$ is a parameter ideal in $T$, using induction hypothesis, it follows that $\ov{e}_k(IT) \geq 0$ and if $\ov{e}_k(IT)=0$, then $\ov{I^{n+k-1}T} \subseteq I^nT$ for all $n.$ Use \cite[Corollary 8]{itoh92} to complete the proof.
\end{proof}

By the Proposition \ref{HSpoly}(b), it is clear that to formulate an analogue of Itoh's theorem for $\ov{e}_k(I)$, $k\leq d$ and $d \geq 3$, we need to relate the vanishing of $\ov{e}_k(I)$ to the equality $\ov{e}_1(I)=\sum_{j=0}^{k-2}\ell(\ov{I^{j+1}}/I\ov{I^j})$.
The following theorem gives a generalization of \cite[Proposition 3.1, 3.2]{cmr}.
%\textcolor{red}{
%\begin{Lemma} \label{polynomial}
%	Let $f(x) \in \mathbb{Z}[x]$ be an integer valued polynomial of positive degree. If there exists an integer $k \gg 0$, such that $f(k)>0,$ then the leading coefficient of $f(x)$ is also positive. 
%\end{Lemma}
%\begin{proof}
%	Since $f(x)$ has a non-zero leading coefficient, suppose the leading coefficient of $f(x)$ is negative. Then for all $n \gg 0,$ we get $f(n)<0.$ This contradicts the given assumption. Thus, the leading coefficient of $f(x)$ is positive. 
%\end{proof}
%}

\begin{Theorem}\label{thm4}
	Let $(R,\m)$ be a $d$-dimensional analytically unramified, Cohen-Macaulay local ring with $d\geq 3$, and $I$ be a parameter ideal such that $\overline{I}=\m$. For some $k \geq 2$, suppose that $\HH^i_{\R(I)_+}(\ov{\R'}(I))_j=0$ for $3\leq i+j\leq k-1$ and $3\leq i\leq d-1$. Let $\ov{e}_k(I)=0$.
	\begin{enumerate}[$(a)$]
		\item Suppose further that $\f=\{\ov{I^n}\}$ satisfies \emph{the condition $HI_{k-2}$}, then 
		\[ \ell\left( \frac{\ov{I^{n+1}}}{I^{n-(k-3)}\ov{I^{k-2}}} \right)\leq t(R){n-(k-3)+d-2\choose d-1} \]
		for all $n \geq k-3$, where $t(R)$ denotes the type of $R$. In particular, $\ell\left(\dfrac{\ov{I^{k-1}}}{I\ov{I^{k-2}}}\right)\leq t(R).$
		\item Suppose $\f$ satisfies $\cc_{r}$ for $2 \leq r \leq k-2$, then 
		\begin{align*}
			\sum_{j=0}^{k-2}\ell\left(\frac{\ov{I^{j+1}}}{I\ov{I^j}}\right)\leq \ov{e}_1(I)\leq \sum_{j=0}^{k-3}\ell\left(\frac{\ov{I^{j+1}}}{I\ov{I^j}}\right)+t(R).
		\end{align*}
		%In particular, if $t(R)\neq \ell\left(\dfrac{\ov{I^{k-1}}}{I\ov{I^{k-2}}}\right)$, then $\displaystyle \ov{e}_1(I)< \sum_{j=0}^{k-3} \ell\left(\frac{\ov{I^{j+1}}}{I\ov{I^j}}\right)+t(R).$
	\end{enumerate}
\end{Theorem}

\begin{proof}
	$(a)$: As $\ov{e}_k(I)=0$, Theorem \ref{vanishingOfe_k} implies that $\ov{I^{n+k-1}}\subseteq I^n$ for all $n \geq 0$, or $\ov{I^{n+2}}\subseteq I^{n-(k-3)}$ for all $n \geq k-3.$ Thus for all $n \geq k-3,$
	\begin{align*}
		\m \ov{I^{n+1}}=\ov{I} \cdot \ov{I^{n+1}}\subseteq \ov{I^{n+2}}\subseteq I^{n-(k-3)}
	\end{align*}
	where the first equality is due to the assumption $\ov{I}=\m$. This implies that
	\begin{align}\label{thm4eq1}
		\ov{I^{n+1}}\subseteq I^{n-(k-3)}:\m.
	\end{align}
	Therefore, for all $n \geq k-3,$
	\begin{align*}
		\ell\left( \frac{\ov{I^{n+1}}}{I^{n-(k-3)}\ov{I^{k-2}}} \right)
		=\ell\left(\frac{\ov{I^{n+1}}}{\ov{I^{n+1}}\cap I^{n-(k-3)}}\right)
		=\ell\left(\frac{\ov{I^{n+1}}+I^{n-(k-3)} }{I^{n-(k-3)}} \right)
		\leq \ell\left(\frac{I^{n-(k-3)}:\m}{I^{n-(k-3)}} \right)
	\end{align*}
	where the first equality is due to the assumption that $\f$ satisfies \emph{the condition $HI_{k-2}$} and the latter inequality is due to equation \eqref{thm4eq1}.
		
	Now $\ell(I^{n-(k-3)}:\m/I^{n-(k-3)})$ is the dimension of the socle of the ring $R/I^{n-(k-3)}.$ In other words, it is the number of irreducible components of $I^{n-(k-3)}$ in $R.$ Write $I=(x_1,\ldots,x_d).$ Using the decomposition of $I^{n-(k-3)}$ as 
	\[ I^{n-(k-3)} = \bigcap_{\substack{c_1+\cdots+c_d=n-(k-3)+d-1 \\ c_1,\ldots,c_d \geq 1}} (x_1^{c_1},\ldots,x_d^{c_d}), \]
	it follows that $\displaystyle \ell(I^{n-(k-3)} : \m / I^{n-(k-3)}) = \alpha \binom{n-(k-3)+d-2}{d-1}$, where $\alpha$ is equal to the number of irreducible components of the ideal $(x_1^{c_1}, \ldots, x_d^{c_d})$ in $R.$ Observe that 
	\[ \alpha= \dim \Soc(R/(x_1^{c_1},\ldots,x_d^{c_d})) = t(R). \] 
	Hence, for all $n \geq k-3,$

	\begin{align*}
		\ell\left(\frac{\ov{I^{n+1}}}{I^{n-(k-3)}\ov{I^{k-2}}} \right)\leq \ell\left(\frac{I^{n-(k-3)}:\m}{I^{n-(k-3)}} \right)=t(R){n-(k-3)+d-2\choose d-1}.
	\end{align*}
	Substituting $n=k-2$ in the above equation gives $\ell\left(\frac{\ov{I^{k-1}}}{I\ov{I^{k-2}}}\right)\leq t(R)$.
		
	$(b)$: Using Theorem \ref{integralclosure}, it follows that $\mathcal{F}$ satisfies \emph{the condition $HI_p$} for all $p \leq k-2.$ Now following the proof as in Proposition \ref{HSpoly}, one can write
	\begin{multline*}
		\ell\left(\frac{R}{\ov{I^{n+1}}}\right)
		= \ell\left(\frac{R}{I}\right)\binom{n+d}{d}-\ell\left(\frac{\ov{I}}{I}\right)  \binom{n+d-1}{d-1} -\\
		\sum_{i=1}^{k-3} \ell\left(\frac{\ov{I^{i+1}}}{I\ov{I^i}}\right) \left[\sum_{j=1}^{i+1} (-1)^{j-1} \binom{i}{j-1}\binom{n+d-j}{d-j} \right] - \ell\left(\frac{\ov{I^{n+1}}}{I^{n-(k-3)}\ov{I^{k-2}}}\right).
	\end{multline*}
	Using part $(a)$, it follows that 
	\begin{multline}\label{thm4eq2}
		\ell\left(\frac{R}{\ov{I^{n+1}}}\right)
		\geq \ell\left(\frac{R}{I}\right)\binom{n+d}{d}-\ell\left(\frac{\ov{I}}{I}\right)  \binom{n+d-1}{d-1} -\\
		\sum_{i=1}^{k-3} \ell\left(\frac{\ov{I^{i+1}}}{I\ov{I^i}}\right) \left[\sum_{j=1}^{i+1} (-1)^{j-1} \binom{i}{j-1}\binom{n+d-j}{d-j} \right] -  t(R){n-(k-3)+d-2\choose d-1}.
	\end{multline}
	From  \cite[Lemma 22]{GMV}, we have
	\begin{align*}
		{n-(k-3)+d-2\choose d-1}={n+d-(k-2)-1\choose d-1}=\sum_{j=1}^{k-1} (-1)^{j-1} \binom{k-2}{j-1}\binom{n+d-j}{d-j}.
	\end{align*}
	Using the above equality in \eqref{thm4eq2}, we have
	\begin{multline*}
		\ell\left(\frac{R}{\ov{I^{n+1}}}\right)
		\geq \ell\left(\frac{R}{I}\right)\binom{n+d}{d}-\ell\left(\frac{\ov{I}}{I}\right)  \binom{n+d-1}{d-1} -\\
		\sum_{i=1}^{k-3} \ell\left(\frac{\ov{I^{i+1}}}{I\ov{I^i}}\right) \left[\sum_{j=1}^{i+1} (-1)^{j-1} \binom{i}{j-1}\binom{n+d-j}{d-j} \right] - t(R)\sum_{j=1}^{k-1} (-1)^{j-1} \binom{k-2}{j-1}\binom{n+d-j}{d-j}.
	\end{multline*}
	This can be simplified as 
	\begin{multline}\label{thm4eq3}
		\ell\left(\frac{R}{\ov{I^{n+1}}} \right)\geq  \ell\left(\frac{R}{I}\right)\binom{n+d}{d} - \beta_1(\f)\binom{n+d-1}{d-1} + \cdots \\
		+ (-1)^{k-2} \beta_{k-2}(\f)\binom{n+d-(k-2)}{d-(k-2)}
		+ t(R)(-1)^{k-1} \binom{n+d-(k-1)}{d-(k-1)}
	\end{multline}
	where $\displaystyle \beta_i(\f) = \sum_{j=i-1}^{k-3} \left[ \binom{j}{i-1}\ell(\ov{I^{j+1}}/I\ov{I^j}) \right] + t(R)\binom{k-2}{i-1}$ for all $i=1,\ldots,k-2.$ Recall that for all $n \gg 0,$ 
	\begin{align} \label{thm4eq4}
	\ell\left(\frac{R}{\ov{I^{n+1}}}\right) = e_0(I) \binom{n+d}{d} - \ov{e}_1(I) \binom{n+d-1}{d-1}+ \cdots + (-1)^d \ov{e}_d(I).
	\end{align}
	Comparing the equations \eqref{poly}, \eqref{thm4eq3} and \eqref{thm4eq4}, we obtain	
	%Now comparing the above inequality with the equation \eqref{poly}, we have
	\begin{align*}
		\sum_{j=0}^{k-2}\ell\left(\frac{\ov{I^{j+1}}}{I\ov{I^j}}\right)
		=\alpha_1(\f)\leq \ov{e}_1(I)\leq \beta_1(\f)
		=\sum_{j=0}^{k-3}\ell\left(\frac{\ov{I^{j+1}}}{I\ov{I^j}}\right)+t(R).
	\end{align*}
\end{proof}

We are now ready to prove the main result of this section.

\begin{Theorem} \label{itohsgeneralization}
	Let $(R,\m)$ be a $d$-dimensional analytically unramified, Cohen-Macaulay local ring with $d\geq 3$, and $I$ be a parameter ideal such that $\overline{I}=\m$. Suppose that for some $2 \leq k \leq d,$ $\f$ satisfies $\cc_{r}$ for $2 \leq r \leq k-2$ and $\ell(\ov{I^{k-1}}/I\ov{I^{k-2}})\geq t(R).$ Then $\ov{e}_k(I)=0$ if and only if $\ov{r}(I)\leq k-1.$ In this case, $\ov{G}(I)$ is Cohen-Macaulay.
\end{Theorem}

\begin{proof}
	Suppose $\ov{e}_k(I)=0$, then using Proposition \ref{thm4}(a) and our assumption, it follows that 
	\begin{align*}
		\ell\left(\frac{\ov{I^{k-1}}}{I\ov{I^{k-2}}}\right)\leq t(R)\leq \ell\left(\frac{\ov{I^{k-1}}}{I\ov{I^{k-2}}}\right).
	\end{align*}
	Thus $t(R)=\ell\left(\dfrac{\ov{I^{k-1}}}{I\ov{I^{k-2}}}\right)$ and hence $\displaystyle \ov{e}_1(R)=\sum_{j=0}^{k-2}\ell\left(\frac{\ov{I^{j+1}}}{I\ov{I^j}}\right)$. It follows that $\ov{r}(I)\leq k-1$ and $\ov{G}(I)$ is Cohen-Macaulay by Proposition \ref{HSpoly}. \\
%	If $t(R)\neq \ell\left(\dfrac{\ov{I^{k-1}}}{I\ov{I^{k-2}}}\right)$, then using Proposition \ref{thm4} we have
%	\begin{align*}
%		\sum_{j=0}^{k-2} \ell\left(\frac{\ov{I^{j+1}}}{I\ov{I^j}}\right)\leq \ov{e}_1(I)
%		< \sum_{j=0}^{k-3} \ell\left(\frac{\ov{I^{j+1}}}{I\ov{I^j}}\right) + t(R)
%		= \sum_{j=0}^{k-3}\ell\left(\frac{\ov{I^{j+1}}}{I\ov{I^j}}\right) + \ell\left(\frac{\ov{I^{k-1}}}{I\ov{I^{k-2}}}\right) + 1.
%	\end{align*}
%	Again the equality $\displaystyle \sum_{j=0}^{k-2}\ell\left(\frac{\ov{I^{j+1}}}{I\ov{I^j}}\right)= \ov{e}_1(I)$ holds and we have $\ov{r}(I)\leq k-1$ and $\ov{G}(I)$ is Cohen-Macaulay by Proposition \ref{HSpoly}.

	The converse is a direct consequence of Proposition \ref{HSpoly}(a).
\end{proof}
\begin{Corollary}[{\rm{\cite[Theorem 3.3]{cmr}}}]
	Let $(R,\m)$ be a $d$-dimensional analytically unramified, Cohen-Macaulay local ring with $d\geq 3$, and $I$ be a parameter ideal such that $\overline{I}=\m$. Suppose that $\ell(\ov{I^{2}}/I\ov{I})\geq t(R).$ Then $\ov{e}_3(I)=0$ if and only if $\ov{r}(I)\leq 2.$ In this case, $\ov{G}(I)$ is Cohen-Macaulay.
\end{Corollary}
\begin{proof}
This is an immediate consequence of the previous theorem with $k=3$.
\end{proof}
The above corollary shows that Theorem \ref{itohsgeneralization}  generalizes parts  of \cite[Theorem 3.3]{cmr} for higher Hilbert coefficients. This is also natural since the techniques presented in this section extend the scope of the techniques in \cite{cmr} and \cite{itoh92}.

\section{Examples}

We illustrate the results in the previous section using the following examples.

\begin{Example}{\rm
	Consider the ring $R=k[[X_0,X_1,\ldots,X_d]]/(X_0^n+\cdots+X_d^n)$, where $\chr k=0$, $d \geq 3$ and $n \leq d.$ Then $R$ is a $d$-dimensional analytically unramified Gorenstein local ring. Let $\m=(X_0,\ldots,X_d)$ be the maximal homogeneous ideal of $R.$ Then the associated graded ring $G(\m) = k[X_0,\ldots,X_d]/(X_0^n+\cdots+X_d^n)$ is a $d$-dimensional Cohen-Macaulay domain. Hence, $\ov{\m^n} = \m^n$ for all $n.$ 
		
	Let $I$ be a parameter ideal of $R$ such that $\ov{I}=\m.$ Then $\f =\{\ov{I^n}\}_{n \in \Z}=\{\m^n\}_{n \in \Z}.$ In this case, the associated graded ring of the filtration, $\ov{G}(I) = G(\m).$ So
	\[ H(\ov{G}(I),t) = H(G(\m),t) = \frac{(1-t^n)}{(1-t)^{d+1}} = \frac{1+t+t^2+\cdots+t^{n-1}}{(1-t)^d}. \]
	This implies that  the postulation number $\ov{n}(I) = n(\m)=n-1-d.$ Since $G(\m)$ is Cohen-Macaulay, we get that the reduction number $\ov{r}(I)=r(\m)=n(\m)+d = n-1$, using \cite[Corollary 3.8]{marleyThesis}. Note that $\ov{e}_n(I) = e_n(\m)=0.$ Moreover, $\ov{e}_k=0$, for all $n \leq k \leq d.$
		
	As $\ov{r}(I) <d$, using \cite[Theorem 2.3]{viet}, it follows that the normal Rees ring $\ov{\R}(I)=\R(\m)$ is Cohen-Macaulay. Hence $\HH^i_J(\ov{\R}(I)) = \HH^i_{(t^{-1},\m t)}(\R(\m))$ is non-zero if and only if $i=d+1.$ Therefore, the filtration $\f$ satisfies \emph{the condition $HI_r$} for all $r.$
}\end{Example}

\begin{Example}{\rm
	Consider the simplicial complex $\Delta$, on 8 vertices, defined by the facets
	\begin{align*}
		\{ &\{5,6,7,8\}, \{2,5,7,8\}, \{1,2,7,8\}, \{3,5,6,8\}, \{2,3,5,8\}, \{1,2,3,8\},\\ &\{4,5,6,7\}, \{2,4,5,7\}, \{1,2,4,7\}, \{3,4,5,6\}, \{2,3,4,5\}, \{1,2,3,4\} \}.
	\end{align*}
	Note that $\Delta$ is a pure simplicial complex and the order in which the facets are written above gives a shelling of $\Delta.$ Hence $\Delta$ is a shellable simplicial complex.
		
	Let $R$ be the Stanley-Reisner ring of $\Delta$. Then $R$ is a 4-dimensional Cohen-Macaulay ring using \cite[Theorem 8.2.6]{herzogHibi}. By localizing at the maximal homogeneous ideal $\m$, we may assume that $R$ is local. Observe that $G(\m) = S/I_{\Delta}$ and  is reduced. This implies that $\ov{\m^n} = \m^n$, for all $n.$ Let $I$ be a parameter ideal in $R$ such that $\ov{I}=\m.$ Then $\f=\{\ov{I^n}\}_{n \in \Z} = \{\m^n\}_{n \in \Z}.$ 
		
	The $f$-vector of the simplicial complex $\Delta$ is $f(\Delta)= (1,8,23,28,12).$ Then the $h$-vector, $h(\Delta)=(1,4,5,2,0).$ Therefore, the Hilbert series
	\[ H(\ov{G}(I),t) = H(G(\m),t) = \frac{1+4t+5t^2+2t^3}{(1-t)^4}. \] 
	This implies that $\ov{e}_4(I) = e_4(\m)=0$ and the postulation number $\ov{n}(I) = n(\m)=-1.$ Since $G(\m)$ is Cohen-Macaulay, using \cite[Corollary 3.8]{marleyThesis}, we get $\ov{r}(I) = r(\m) = 3.$ Further, this implies that the normal Rees ring $\ov{\R}(I) = \R(\m)$ is Cohen-Macaulay and so the filtration $\f$ satisfies \emph{the condition $HI_r$} for all $r.$
}\end{Example} 
\section{The condition $HI_r$ for any $I$-admissible filtration}

In this section, we discuss the necessary conditions for any $I$-admissible filtration $\f$ to satisfy \emph{the condition $HI_r$}. Though an exact analogue of Theorem \ref{integralclosure} is sought after, extra conditions are required.

\begin{Theorem} \label{HIforfiltration}
	Let $R$ be a $d$-dimensional Cohen-Macaulay local ring and let $I=(a_1,\dots,a_g)$ be an $R$-ideal generated by a regular sequence of height $g \ge 1.$ Let $\f=\{ I_n \}$ be an $I$-admissible filtration. If $\f$ satisfies $\cc_{r}$ and $\depth G(\f) \ge g-1$, then $\f$ satisfies \emph{the condition $HI_r$}.
\end{Theorem}

\begin{proof}
	We proceed by induction on $g.$ Let $g=1$ and $I=(a).$ Let $u=t^{-1}$ and consider the \v{C}ech complex
	\[0 \longrightarrow \R'(\f) \overset{f}\longrightarrow \R'(\f)_{at} \times \R'(\f)_u \overset{g}\longrightarrow \R'(\f)_{atu} \longrightarrow 0\]
	where 
	\[f(e)=\left(\frac{e}{1},\frac{e}{1}\right) \text{ \ and \ } g\left(\frac{b}{(at)^{n_1}} , \frac{c}{u^{n_2}}\right) = \frac{b}{(at)^{n_1}}-\frac{c}{u^{n_2}}.\]
	We show that $I_{n+r}\cap I^{n}\subseteq I_{r}I^{n}$, since the other inclusion is clear. Let $z\in I_{n+r}\cap (a)^{n}.$ Write $z=ba^{n}$, where $b\in R.$ Consider 
	\[\left(\frac{zt^{n+r}}{(at)^{n}} , \frac{b}{u^r}\right) \in [\R'(\f)_{at} \times \R'(\f)_{u}]_r.\]
	Observe that 
	\[g\left(\frac{zt^{n+r}}{(at)^{n}} , \frac{b}{u^r}\right) = \frac{zt^{n+r}}{(at)^{n}} - \frac{b}{u^r} = \frac{zt^{n}-ba^{n}t^{n}}{(at)^{n}u^r}=0.\]
	This implies that 
	\[\left(\frac{zt^{n+r}}{(at)^{n}} , \frac{b}{u^r}\right) \in [\ker(g)]_r.\]
	Since $\HH^1_J(\R'(\f))_r=0$, $\ker(g)_r = \im(f)_r.$ It follows that there exists $e \in I_r$ such that 
	\[ f(et^r)=\left(\frac{et^r}{1},\frac{et^r}{1}\right)=\left(\frac{zt^{n+r}}{(at)^{n}} , \frac{b}{u^r}\right). \]
	Therefore there exists $s \in \N$ such that $(at)^s(ea^{n}t^{n+r} - zt^{n+r})=0$ in $\R'(\f).$ As $a$ is a nonzerodivisor in $R$, we get $z=ea^{n} \in I_rI^{n}.$ Hence, $I_{n+r}\cap I^{n} = I_{r}I^{n}$, for all $n \ge 0.$
	
	Let $g \ge 2$ and let $I=(a_1,\dots,a_g)$, where $a_1,\dots,a_g$ is an $R$-regular sequence. Then $R/(a_1)$ is Cohen-Macaulay and $I/(a_1)$ is ideal generated by $R/(a_1)$-regular sequence such that $\hgt(I/(a_1)) \ge 1.$ Also, $\f/(a_1) = \{ (I_n+(a_1))/(a_1) \}$ is an $I/(a_1)$-admissible filtration and $\depth G(\f/(a_1t)) \ge g-2.$ Consider the Rees ring of the filtration $\f/(a_1)$, given by 
	\[\R(\f/(a_1)) = \bigoplus_{n\in \Z} \left(\frac{I_n+(a_1)}{(a_1)} \right) t^n 
	= \bigoplus_{n\in \Z} \left(\frac{I_n}{(a_1) \cap I_n} \right) t^n 
	= \bigoplus_{n\in \Z} \left(\frac{I_n}{a_1I_{n-1}} \right) t^n \simeq \frac{\R'(\f)}{(a_1t)}.\]
	The latter equality is true because the image of $a_1$ in $G(\f)$ is a nonzerodivisor. Consider the exact sequence
	\[ 0 \rightarrow \R'(\f)(-1) \overset{a_1t}\longrightarrow \R'(\f) \longrightarrow \frac{\R'(\f)}{(a_1t)} \rightarrow 0. \]
	This gives a long exact sequence of local cohomology modules
	\[ \cdots \rightarrow \HH^i_J(\R'(\f))_j \rightarrow \HH^i_J(\R'(\f)/(a_1t))_j \rightarrow \HH^{i+1}_J(\R'(\f))_{j-1} \rightarrow \cdots \]
	Using the change of ring principle, 
	\[ \cdots \rightarrow \HH^i_J(\R'(\f))_j \rightarrow \HH^i_{J/(a_1t)}(\R'(\f)/(a_1t))_j \rightarrow \HH^{i+1}_J(\R'(\f))_{j-1} \rightarrow \cdots \]
	As $\HH^i_J(\R'(\f))_{j}=0$ for all $i,j$ such that $i+j=r+1$ and $1 \le i \le g$, we get $\HH^i_{J/(a_1t)}(\R'(\f)/(a_1t))_{j}=0$ for all $i,j$ such that $i+j=r+1$ and $1 \le i \le g-1=\hgt(I/(a_1)).$ Using induction hypothesis, it follows that for all $n\geq 0$, 
	\begin{equation} \label{1}
	(a_1)+I_{n+r}\cap I^{n} = I_{r}I^{n} + (a_1).
	\end{equation}
	We show that $I_{n+r}\cap I^{n}\subseteq I_{r}I^{n}$ for all $n\geq 0$, as the other inclusion is clear. Let $z \in I_{n+r}\cap I^{n}.$ Using equation (\ref{1}), write $z=a+ba_1$, where $a \in I_rI^{n}$ and $b \in R.$ Hence, $b \in (z-a) \colon a_1.$ In particular, as the image of $a_1$ in $G(\f)$ is a nonzerodivisor, it follows that 
	\[b \in (I_{n+r} \cap I^{n}) \colon a_1 = (I_{n+r} \colon a_1) \cap (I^{n} \colon a_1) = I_{n+r-1} \cap I^{n-1}.\]
	This implies that $z \in I_rI^{n} + (I_{n+r-1} \cap I^{n-1})a_1.$ Induction on $n$ gives $I_{n+r-1} \cap I^{n-1} = I^{n-1}I_r.$ Hence, $z \in I_rI^{n}.$ 
\end{proof}

\begin{Remarks}{\rm
		(1) Comparing Theorem \ref{HIforfiltration} and Theorem \ref{integralclosure}, observe that, for the normal filtration, we do not require the assumption on the depth of the associated graded ring $G(\f)$. \\
		(2) Let $g = \hgt(I).$ Since we assume $\depth G(\f) \geq g-1$, it follows that $\depth \R'(\f) \geq g.$ Hence $\HH^i_J(\R'(\f))=0$ for all $i \leq g-1.$ Thus, if $\HH^g_J(\R'(\f))_{r+1-g}=0$ and $\depth G(\f) \geq g-1$, then $\f$ satisfies {\em the condition $HI_r.$}
}\end{Remarks}

\bibliographystyle{plain}  		
\bibliography{gmv2019}

\begin{thebibliography}{10}

\bibitem{brunsHerzog}
Winfried {B}runs and J\"urgen Herzog.
\newblock {\em {C}ohen-{M}acaulay rings}, volume~39 of {\em Cambridge Studies
  in Advanced Mathematics}.
\newblock Cambridge University Press, Cambridge, 1993.

\bibitem{cmr}
Alberto Corso, Claudia Polini, and Maria~Evelina Rossi.
\newblock Bounds on the normal {H}ilbert coefficients.
\newblock {\em Proc. Amer. Math. Soc.}, 144(5):1919--1930, 2016.

\bibitem{GMV}
Kriti Goel, Vivek Mukundan, and J.~K. Verma.
\newblock Tight closure of powers of ideals and tight {H}ilbert polynomials.
\newblock {\em arXiv preprint arXiv:1806.07522}, pages 1--24, 2018.

\bibitem{herzogHibi}
J\"{u}rgen Herzog and Takayuki Hibi.
\newblock {\em Monomial ideals}, volume 260 of {\em Graduate Texts in
  Mathematics}.
\newblock Springer-Verlag London, Ltd., London, 2011.

\bibitem{hochsterGeneralGradeReductions}
M.~Hochster.
\newblock Properties of {N}oetherian rings stable under general grade
  reduction.
\newblock {\em Arch. Math. (Basel)}, 24:393--396, 1973.

\bibitem{hongUlrich}
Jooyoun Hong and Bernd Ulrich.
\newblock Specialization and integral closure.
\newblock {\em J. Lond. Math. Soc. (2)}, 90(3):861--878, 2014.

\bibitem{Huckaba1}
Sam Huckaba.
\newblock A {$d$}-dimensional extension of a lemma of {H}uneke's and formulas
  for the {H}ilbert coefficients.
\newblock {\em Proc. Amer. Math. Soc.}, 124(5):1393--1401, 1996.

\bibitem{HunekeHuckaba99}
Sam Huckaba and Craig Huneke.
\newblock Normal ideals in regular rings.
\newblock {\em J. Reine Angew. Math.}, 510:63--82, 1999.

\bibitem{huckabaMarley}
Sam Huckaba and Thomas Marley.
\newblock Hilbert coefficients and the depths of associated graded rings.
\newblock {\em Journal of the London Mathematical Society}, 56(1):64--76, 1997.

\bibitem{huneke1987}
Craig Huneke.
\newblock Hilbert functions and symbolic powers.
\newblock {\em Michigan Math. J.}, 34(2):293--318, 1987.

\bibitem{itoh88}
Shiroh Itoh.
\newblock Integral closures of ideals generated by regular sequences.
\newblock {\em J. Algebra}, 117(2):390--401, 1988.

\bibitem{itoh92}
Shiroh Itoh.
\newblock Coefficients of normal {H}ilbert polynomials.
\newblock {\em J. Algebra}, 150(1):101--117, 1992.

\bibitem{kuma}
Manoj Kummini and Shreedevi Masuti.
\newblock On conjectures of {I}toh and of {L}ipman on the cohomology of
  normalized blow ups.
\newblock {\em arxiv preprint arXiv:1507.03343}, pages 1--17, 2015.

\bibitem{mandal2012chern}
Mousumi Mandal and J.~K. Verma.
\newblock On the {C}hern number of {$I$}-admissible filtrations of ideals.
\newblock {\em J. Commut. Algebra}, 4(4):577--589, 2012.

\bibitem{marley89}
Thomas Marley.
\newblock The coefficients of the {H}ilbert polynomial and the reduction number
  of an ideal.
\newblock {\em Journal of the London Mathematical Society}, 2(1):1--8, 1989.

\bibitem{marleyThesis}
Thomas~John Marley.
\newblock {\em Hilbert functions of ideals in {C}ohen-{M}acaulay rings}.
\newblock ProQuest LLC, Ann Arbor, MI, 1989.
\newblock Thesis (Ph.D.)--Purdue University.

\bibitem{ooishi}
Akira Ooishi.
\newblock Castelnuovo's regularity of graded rings and modules.
\newblock {\em Hiroshima Math. J.}, 12(3):627--644, 1982.

\bibitem{reesAUrings}
D~Rees.
\newblock A note on analytically unramified local rings.
\newblock {\em Journal of the London Mathematical Society}, 1(1):24--28, 1961.

\bibitem{viet}
Duong~Qu\^oc Vi\^et.
\newblock A note on the {C}ohen-{M}acaulayness of {R}ees algebras of
  filtrations.
\newblock {\em Comm. Algebra}, 21(1):221--229, 1993.

\end{thebibliography}
\end{document}